\documentclass[12pt, draft]{article}
\usepackage{amsmath}
\usepackage{amssymb}
\usepackage{latexsym}
\allowdisplaybreaks[4]
\renewcommand{\theequation}{\arabic{section}.\arabic{equation}}
\newtheorem{Thm}{Theorem}[section]
\newtheorem{Prop}[Thm]{Proposition}
\newtheorem{Lemma}[Thm]{Lemma}

\newtheorem{Remark}[Thm]{Remark}
\def\Bbb R{{\rm \bf R}}
\renewcommand{\v}[1]{\vert #1 \vert}
\newcommand{\V}[1]{\Vert #1 \Vert}

\title{On finding a buried obstacle in a layered medium
via the time domain enclosure method}
\author{Masaru IKEHATA\thanks{
Laboratory of Mathematics,
Graduate School of Engineering,
Hiroshima University, Higashihiroshima 739-8527, JAPAN;
ikehata@hiroshima-u.ac.jp
}
\and 
Mishio KAWASHITA\thanks{Department of Mathematics,
Graduate School of Sciences,
Hiroshima University, Higashihiroshima 739-8526, JAPAN;
kawasita@hiroshima-u.ac.jp}
}
\date{}
\begin{document}
\maketitle
\begin{abstract}
An inverse obstacle problem for the wave equation in a two layered medium 
is considered.  It is assumed that the unknown obstacle is penetrable and 
embedded in the lower half-space.
The wave as a solution of the wave equation is generated by an initial data
whose support is in the upper half-space and observed at the same place as the support
over a finite time interval.  From the observed wave an indicator function 
in the time domain enclosure method is constructed.
It is shown that, one can find some information about 
the geometry of the obstacle together with the qualitative property in the asymptotic behavior of the indicator function.
\end{abstract}


\vskip 1truepc
\setcounter{section}{0}
\section{Introduction}\label{Introduction}
\vskip-6pt\noindent

The problem of finding an obstacle embedded or hidden in a complicated environment
by using the electromagnetic wave appears in, for example,
ground penetrating or subsurface radar \cite{DGS} and
the through-wall imaging \cite{BA}. 

In this paper, we consider such type of the problems in a simplest, however, important mathematical model which employs
a wave governed by a scalar wave equation in a two homogeneous layered medium over a finite time interval.

Let $\Bbb R^3_{\pm}=\{x=(x_1,x_2,x_3)\in\Bbb R^3\,\vert\,\pm x_3>0\}$.
Consider $\gamma_0\in L^{\infty}(\Bbb R^3)$ given by 
$$\displaystyle
\gamma_0(x)=\left\{
\begin{array}{lr}
\displaystyle \gamma_+, & \quad\text{if $x_3>0$,}
\\
\displaystyle \gamma_-, & \quad\text{if $x_3<0$,}
\end{array}
\right.
$$
where $\gamma_{\pm}$ are positive constants.

Let $D$ be a bounded open subset of $\Bbb R^3_{-}$ with a $C^2$-boundary and satisfy
$\overline D\subset\Bbb R^3_{-}$.  
Consider $\gamma\in L^{\infty}(\Bbb R^3)$ given by
$$\displaystyle
\gamma(x)=\left\{
\begin{array}{ll}
\displaystyle \gamma_0(x)I_3, & \quad\text{if $x\in\Bbb R^3\setminus D$,}
\\
\displaystyle \gamma_0(x)I_3+h(x), & \quad\text{if $x\in D$, }
\end{array}
\right.
$$
where $h=h(x)$, $x\in D$ is a real symmetric $3\times 3$-matrix valued function
and satisfies that: all the components of $h$ are essentially bounded on $D$
; there exists a positive constant $C$ such that
$(\gamma_0(x)I_3+h(x))\xi\cdot\xi\ge C\vert\xi\vert^2$ for all $\xi\in\Bbb R^3$ and a.e. $x\in D$.

Let $0<T<\infty$.
Given $f\in L^2(\Bbb R^3)$ let $u=u_f(x,t)$ be the weak solution of the initial value problem:
\begin{equation}
\left\{
\begin{array}{ll}
\displaystyle
(\partial_t^2-\nabla\cdot\gamma\nabla) u=0 & \text{in}\,\Bbb R^3\times\,]0,\,T[,
\\
\displaystyle
u(x,0)=0 & \text{on}\,\Bbb R^3,
\\
\displaystyle
\partial_tu(x,0)=f(x) & \text{on}\,\Bbb R^3.
\end{array}
\right.
\label{1.1}
\end{equation}
Note that the solution class is taken from \cite{DL}.  See also \cite{IE} for its detailed description.

We consider the following problem:

\vskip1pc

{\bf\noindent Problem.}
Assume that $\gamma_{+}\not=\gamma_{-}$.
Fix a large $T$ (to be determined later).
Assume that $\gamma_0$ is {\it known} and that both $D$ and $h$ are {\it unknown}.
Let $B$ be an open ball whose closure is contained in $\Bbb R_{+}^3$.
Fix a $f\in L^2(\Bbb R^3)$ with $\text{supp}\,f\subset\overline B$ and 
satisfying that there exists a positive constant $C_0$ such that $f(x)\ge C_0$ a.e. $x\in B$ or
$-f(x)\ge C_0$ a.e. $x\in B$.
Generate $u=u_f$ of the solution of (\ref{1.1}) by the $f$.
Extract information about the location and shape of $D$ from the measured data $u$ on $B$
over the time interval $]0,\,T[$.

\vskip1pc

It should be emphasized that the problem asks us to extract information about unknown obstacle $D$
from a {\it single wave} observed over a {\it finite time interval} at the same place where the wave is generated.
There are some studies which consider the {\it time harmonic reduced case} in a two layered medium.
See \cite{LZ} for uniqueness issue of impenetrable obstacles
using {\it infinitely many} incident fields; \cite{LLLL} a reconstruction scheme of an impenetrable obstacle using a far-field pattern
corresponding to a single incident plane wave; \cite{DEKPS} study of a direct problem with an application to mine detection
and propose a numerical reconstruction scheme using near field measurements corresponding to finitely many incident sources.
Clearly our problem formulation is different from their one and to our best knowledge there is no result for the problem.

In \cite{IE} Ikehata has considered the case when $\gamma_+=\gamma_{-}(=1)$
and the wave is observed on a closed surface $S$ over a finite time interval which encloses 
the obstacle.  He assumed that $\gamma$ satisfies
one of the following two conditions:

\noindent
(A1)  there exists a positive constant $C'$ such that $-h(x)\xi\cdot\xi\ge C'\vert\xi\vert^2$
for all $\xi\in\Bbb R^3$ and a.e. $x\in D$;

\noindent
(A2)  there exists a positive constant $C'$ such that $h(x)\xi\cdot\xi\ge C'\vert\xi\vert^2$
for all $\xi\in\Bbb R^3$ and a.e. $x\in D$.

\noindent
In Theorem 1.2 of \cite{IE} he showed that if $B$ is outside surface $S$, then one can extract the
distance $\text{dist}\,(D, B)=\inf_{x\in D, y\in B}\vert x-y\vert$ from the observed wave
and also can distinguish whether (A1) or (A2) is satisfied by using the signature of an
indicator function computed from the observed wave.  This is the beginning of the multi-dimensional version of the time domain enclosure
method \cite{I4} for inverse obstacle scattering in the time domain.  
In \cite{IEO2} this idea has been extended to the case when
the wave is observed on the same place as the support of an initial data.  This is a version of the near field inverse back-scattering
problem.  One can easily transplant the results in \cite{IE} to this case as pointed out in Subsection 1.3 of \cite{IEO2}.
However, the case when $\gamma_{+}\not=\gamma_{-}$ is not trivial.
Clearly this is a quite interesting case from practical and mathematical point of view.
The unknown obstacle is embedded in the lower half-space which has a different refractive index
from the upper half-space.  Thus the wave generated by an initial data produces reflected and refracted waves at the interface.
The produced refracted wave hits the surface of the obstacle and generates
reflected and refracted waves.  How can one extract information about the geometry of the obstacle from
the observed wave?  The aim of this paper is to extend the previous result to the case when $\gamma_{+}\not=\gamma_{-}$
using the enclosure method in the time domain.  

Now let us describe our main result.
Let $\tau>0$. Let $u$ be the solution of (\ref{1.1}). Define
\begin{align}
w=w_f(x,\tau)
=\int_0^T e^{-\tau T}u(x,t)dt, 
\qquad 
x\in\Bbb R^3.
\label{definition of w}
\end{align}
Let $v\in H^1(\Bbb R^3)$ be the weak solution of
\begin{equation}
\begin{array}{ll}
\displaystyle
(\nabla\cdot\gamma_0\nabla-\tau^2)v+f=0 & 
\displaystyle
\text{in}\,\Bbb R^3.
\end{array}
\label{1.2}
\end{equation}
Define
$$
\displaystyle
I_f(\tau,T)
=\int_{\Bbb R^3}f(w-v)dx.
$$
We call the function $\tau\longmapsto I_f(\tau,T)$ the {\it indicator function}.
Note that this symbol follows from that of \cite{IWALL}.

Define
\begin{equation*}
\displaystyle
l(D,B)=\inf_{x\in D,\,y\in B}l(x, y),
\end{equation*}
where
\begin{equation}
l(x, y) = \inf_{z' \in \Bbb R^2}l_{x, y}(z')
\label{definition of l(x, y)}
\end{equation}
and
\begin{align}
l_{x, y}(z') = \frac{1}{\sqrt{\gamma_-}}\v{\tilde{z}'-x}
+\frac{1}{\sqrt{\gamma_+}}\v{\tilde{z}'-y}
\quad(\tilde{z}' = (z_1, z_2, 0), z' = (z_1, z_2)).
\label{the path length for Snell's law}
\end{align}
The quantity $l(D,B)$ corresponds to the {\it optical distance} or {\it optical path length}  between $B$ and $D$ in optics
and it is easy to see that we have
$$\displaystyle
l(D,B)=l(D,p)-\frac{\eta}{\sqrt{\gamma_+}},
$$
where $p$ and $\eta$ denote the center and radius of $B$, respectively and
$$\displaystyle
l(D,p)=\inf_{x\in D}\,l(x,p).
$$
Thus the unknown obstacle $D$ is contained in the set 
$$
\displaystyle
E(D;B,\gamma_+,\gamma_{-})\equiv\left\{x\in\Bbb R^3_{-}\,\vert\,l(x,p)>l(D,B)+\frac{\eta}{\sqrt{\gamma_+}}\right\}.
$$

The following theorem is the main result of this paper.

\begin{Thm}\label{Theorem 1.1} Assume that $\gamma_{+}<\gamma_{-}$.
Then, we have:
$$
\displaystyle
\lim_{\tau\rightarrow\infty}e^{\tau T}I_{f}(\tau,T)
=
\left\{
\begin{array}{ll}
\displaystyle
0,
& \text{if $T<2l(D,B)$,}\\
\\
\displaystyle
\infty,
&
\text{if $T>2l(D,B)$ and $\gamma$ satisfies (A1),}\\
\\
\displaystyle
-\infty,
&
\text{if $T>2l(D,B)$ and $\gamma$ satisfies (A2).}
\end{array}
\right.
$$
Moreover, if $\gamma$ satisfies (A1) or (A2), then
for all $T>2l(D,B)$
\begin{equation*}
\displaystyle
\displaystyle
\lim_{\tau\rightarrow\infty}\frac{1}{\tau}
\log\left\vert I_{f}(\tau,T) \right\vert
=-2l(D,B).
\end{equation*}
\end{Thm}

From Theorem \ref{Theorem 1.1} we see that the $T$ in the problem should be an arbitrary number satisfying
$T>2l(D,B)$.  We think that this is optimal.
From the indicator function one gets the value $l(D,B)$ and thus the set $E(D;B,\gamma_+,\gamma_{-})$
which {\it encloses} $D$.
Moreover, one can distinguish whether unknown obstacle $D$ satisfies
(A1) or (A2) which is a {\it qualitative property} of $D$ relative to the surrounding background medium, 
by checking the asymptotic behavior of the indicator function.

\begin{Remark}\label{Why gamma_+ < gamma_-}
Intuitively, any signal emanating from $B$ reaches $D$. 
For these signals to go back to the upper side, we need to catch the refracted waves of the reflected waves
by $D$. Hence we need to take measurement in $B$ at least till time $2l(D,B)$ for which the fastest 
signals may come back.
To check whether signals are exactly coming back, we need to take account of total reflection waves.
Assumption $\gamma_+ < \gamma_-$ means that the propagation speed of waves in the upper side is slower
than that of the lower side. Hence, there is no total reflected wave for the incident waves 
from the lower side. This is the case that we do not need to take care of it. 
Mathematically, this is appeared as a difficulty for obtaining asymptotics of the refracted wave.  
As is in (\ref{asymptotics of the refracted wave}) and 
(\ref{asymptotics of the gradient of the refracted wave}) below, it is relatively simple 
since it does not contain waves for total reflection.
\end{Remark}

The proof of Theorem \ref{Theorem 1.1} employs two important facts.
The first one is the following lemma.
\begin{Lemma}\label{Lemma 1.1}
We have, as $\tau\longrightarrow\infty$
\begin{equation}
\displaystyle
I_f(\tau,T)
\ge
\int_{\Bbb R^3}(\gamma_0I_3-\gamma)\nabla v\cdot\nabla v dx+O(\tau^{-1}e^{-\tau T})
\label{1.5}
\end{equation}
and
\begin{equation}
\displaystyle
I_f(\tau,T)
\le \int_{\Bbb R^3}\gamma_0(\gamma_0I_3-\gamma)\gamma^{-1/2}
\nabla v\cdot\gamma^{-1/2}\nabla vdx+O(\tau^{-1}e^{-\tau T}).
\label{1.6}
\end{equation}
\end{Lemma}

For the proof see Appendix.
Combining (\ref{1.5}) and (\ref{1.6}) under the assumption (A1) or (A2), we can easily see that 
Theorem \ref{Theorem 1.1} can be proved if one has the following fact concerning with the asymptotic 
behavior of $\nabla v$ on $D$ as $\tau\longrightarrow\infty$.

\begin{Thm}\label{Theorem 1.2} 
Assume that $\partial D$ is $C^1$ and that $\gamma_{+}<\gamma_{-}$.
Then, there exist positive numbers $C$ and $\tau_0$ such that, for all $\tau\ge\tau_0$
we have
\begin{equation}
\displaystyle
C^{-1}\tau^{-4}e^{-2\tau l(D,B)}
\le\int_D\vert\nabla v(x)\vert^2\,dx\le
C\tau^2e^{-2\tau l(D,B)}.
\label{1.7}
\end{equation}
\end{Thm}

This is the second important fact.
We found that the proof of Theorem \ref{Theorem 1.2} is not a simple matter
and thus the remaining part of this paper is devoted to the proof.
In this sense, the main contribution of this paper to the enclosure method in the time domain
is the establishment of the estimate (\ref{1.7}).  Note that in \cite{KLS, B} one can find 
some formal asymptotic computation of the solution of (\ref{1.2}), however, we do not know 
whether or not their formal theory enables us to derive estimate (\ref{1.7}).

\par

In \cite{IWALL}, Ikehata considered a mathematical model
of the through-wall imaging by using the enclosure method in the time domain.  
Originally the governing equation should be the Maxwell system,
however, as a first step, it is assumed that the governing equation
is given by the single wave equation 
$\alpha(x)\partial_t^2u-\triangle u=0$ in $\Bbb R^3\times\,]0\,,T[$
with the initial data $u(x,0)=0$ and $\partial_t u(x,0)=f(x)$ in $\Bbb R^3$. 
The assumption on the function $\alpha\in L^{\infty}(\Bbb R^3)$ is that: $\alpha$ has a positive lower bound
in $\Bbb R^3$ and takes the form
$$\displaystyle
\alpha(x)=\left\{
\begin{array}{ll}
\displaystyle \alpha_0(x), & \quad\text{if $x\in\Bbb R^3\setminus D$,}
\\
\displaystyle \alpha_0(x)+h(x), & \quad\text{if $x\in D$},
\end{array}
\right.
$$
where the function $\alpha_0$ is essentially bounded in $\Bbb R^3$ with a positive lower bound;
$D$ is an arbitrary bounded open set of the whole space with a Lipschitz boundary; 
the function $h(x)$ or $-h(x)$ has a positive essential infimum on $D$.

Remarkably enough, in \cite{IWALL} a higher regularity more than 
the essential boundedness of $\alpha_0$ is not assumed. 
Thus the model covers various background media such as multilayered media 
with complicated interfaces or unions 
of various domains with different refractive indexes.  He showed that an indicator 
function computed from the wave observed on the same place as the support of 
an initial data yields {\it lower and upper estimates}
of the distance $\text{dist}\,(D,B)$ together with a criterion whether 
$\text{ess.}\inf_{x\in D}h(x)>0$
or $\text{ess.}\inf_{x\in D}(-h(x))>0$ provided $\alpha_0$ is {\it known}. 
The result is based on a system of inequalities similar to (\ref{1.5}) and (\ref{1.6}) in 
Lemma \ref{Lemma 1.1} and explicit upper and lower estimates of the solution of the equation
$$\begin{array}{ll}
\displaystyle
\Delta v-\alpha_0\tau^2 v+\alpha_0f=0 & \text{in}\,\Bbb R^3.
\end{array}
$$

In contrast to this result, Theorem \ref{Theorem 1.1} tells us that
one can extract the {\it exact value}
of the {\it optical} distance $l(D,B)$ from the asymptotic behavior of the indicator function
under the assumption that the background medium consists
of two isotropic {\it homogeneous} layered media.
It would be possible to apply the idea of the derivation of the estimate in Theorem \ref{Theorem 1.2} 
to the case when
$\alpha_0$ takes two different constants, $\alpha_+$ in $x_3>0$ and $\alpha_-$ in $x_3<0$ provided 
$\alpha_+>\alpha_-$ which corresponds to $\gamma_+<\gamma_-$.
However, a typical case to be considered for the Maxwell system is: the upper layer consists of air
and the lower of material, like soil, wall, etc., see \cite{DGS} and \cite{BA}. 
In our problem setting this corresponds to the case when $\gamma_+>\gamma_-$.
Developing an analysis that covers this case together with application to the Maxwell system belongs to our next project.  
See also \cite{IR} for a survey on recent results for inverse obstacle scattering via the time domain enclosure method
and \cite{IK5} for applications to the inverse boundary value problems for the heat equation
in three-dimensional space.

The outline of this paper is as follows.
Let $\Phi_\tau(x, y)$ be the fundamental solution of (\ref{1.2}), which satisfies
\begin{align*}
\nabla_x\cdot(\gamma_0(x)\nabla_x \Phi_\tau(x, y)) -\tau^2\Phi_\tau(x, y)
+\delta(x-y) = 0
\qquad \text{in } \Bbb R^3.
\end{align*}
Since the solution $v$ of (\ref{1.2}) is written by the convolution of $f$ and 
$\Phi_\tau(x, y)$ as
$$
v(x) = \int_{B}\Phi_\tau(x, y)f(y)dy,
$$
we obtain
\begin{align}
\int_{D}\v{\nabla_xv(x)}^2dx
= \int_{B}dy\int_{B}d{\xi}f(y){f(\xi)}
\int_{D}\nabla_x\Phi_\tau(x, y)\cdot{\nabla_x\Phi_\tau(x, \xi)}dx. 
\label{the integral of the gradient of v in D}
\end{align}
Thus, we need to know an asymptotic behavior of $\nabla_x\Phi_\tau(x, y)$
as $\tau \to \infty $ for $x = (x', x_3)$ with $x_3 < 0$, $x' \in \Bbb R^2$
and $y \in B$.
\par

The first step for obtaining the asymptotic behavior of 
(\ref{the integral of the gradient of v in D}) is to show that
the fundamental solution $\Phi_\tau(x, y)$ for $x_3 < 0$ is given by
\begin{align}
\Phi_\tau(x, y)
= \frac{\tau}{4\pi\gamma_+}\int_{\Bbb R^2}
E^{\gamma_-}_{\tau}(x, z')
\frac{e^{-\tau\v{\tilde{z}'- y}/\sqrt{\gamma_+}}}{\v{\tilde{z}'- y}}dz',
\label{expression of the FS for x_3 < 0}
\end{align}
which is given in Section \ref{The refracted part of the fundamental solution}.
In (\ref{expression of the FS for x_3 < 0}), 
$E^{\gamma_-}_{\tau}(x, z')$ is a function 
given by 
\begin{equation}
E^{\gamma_-}_{\tau}(x, z')
= \frac{\tau}{(2\pi)^3}
\int_{\Bbb R^3}{e}^{i\tau\xi\cdot(x-\tilde{z}')}
\frac{1}{\gamma_-\xi^2+1}R_-(\v{\xi'})d{\xi}
\qquad(x_3 < 0), 
\label{the refracted wave}
\end{equation}
where $\tilde{z}' = (z', 0)$ $(z' \in \Bbb R^2)$ is the point on 
the transmission boundary $\partial\Bbb R^3_\pm$ and
$R_-(\v{\xi'})$ is a function of $\v{\xi'}$ standing for 
the transmission coefficient given by 
\begin{align}
R_-(\rho) = 
\frac{4\gamma_+\gamma_-\sqrt{1/\gamma_+ + \rho^2}\sqrt{1/\gamma_- + \rho^2}}
{\gamma_+\sqrt{1/\gamma_+ + \rho^2}+\gamma_-\sqrt{1/\gamma_- + \rho^2}}
\qquad(\rho \geq 0).
\label{the form of R}
\end{align}
Note that $E^{\gamma_-}_{\tau}(x, z')$ can be interpreted
as the refracted part of the fundamental solution.
\par
We put 
$$
E^{\gamma_+, 0}_\tau(x, y) = \frac{1}{4\pi\gamma_+}\frac{e^{-\tau\v{x - y}/\sqrt{\gamma_+}}}{\v{x - y}}
\quad(x \neq y, \tau > 0),
$$ 
which is a fundamental solution for the equation corresponding to the case that
there is no transmission boundary,
i.e. $\gamma_- = \gamma_+$, and given by
\begin{equation}
E^{\gamma_+, 0}_\tau(x, y) = \frac{1}{(2\pi)^3}\int_{\Bbb R^3}
{e}^{i\xi\cdot(x-y)}
\frac{1}{\gamma_+\xi^2+\tau^2}
d{\xi}
= \frac{\tau}{(2\pi)^3}\int_{\Bbb R^3}{e}^{i\tau\xi\cdot(x-y)}
\frac{1}{\gamma_+\xi^2+1}d\xi.
\label{Fourier expression of FS for the free case}
\end{equation}
Thus, (\ref{expression of the FS for x_3 < 0}) stands for 
the refraction phenomena 
by the transmission boundary $\partial\Bbb R^3_\pm$.

\par

As in Proposition \ref{Asymptotics of refracted parts for case 1}, 
for any $N \in {\bf N} $, 
the refracted wave $E^{\gamma_-}_{\tau}(x, z')$
is of the form:
\begin{align}
E^{\gamma_-}_{\tau}(x, z')
= \frac{e^{-{\tau\v{x-\tilde{z}'}/\sqrt{\gamma_-}}}}{4\pi\gamma_-\v{x - \tilde{z}'}}
\Big(\sum_{j = 0}^{N-1}E_{j}(x-\tilde{z}')
\Big(\frac{\sqrt{\gamma_-}}{\tau\v{x - \tilde{z}'}}\Big)^j+\tilde{E}_N(x, z'; \tau)
\Big),
\label{asymptotics of the refracted wave}
\end{align}
where each $ E_{j}$ is a $C^\infty$ function in 
$\Bbb R^3_-$, and $\tilde{E}_N(x, z'; \tau)$ is a continuous for $x \in \Bbb R^3_-$ and $z' \in \Bbb R^2$, 
and satisfies
$$
\v{\tilde{E}_N(x, z'; \tau)} \leq C_N\Big(\frac{\sqrt{\gamma_-}}{\tau\v{x - \tilde{z}'}}\Big)^N
\qquad(x \in \Bbb R^3_-, z' \in \Bbb R^2, \tau > 1).
$$
For the gradient $\nabla_xE^{\gamma_-}_{\tau}(x, z')$, 
as is in Proposition \ref{Asymptotics of refracted parts for case 1} and 
Remark \ref{Asymptotics of gradient of refracted parts for case 1}, we have 
\begin{align}
\nabla_xE^{\gamma_-}_{\tau}(x, z')
= \frac{-{\tau}e^{-{\tau\v{x-\tilde{z}'}/\sqrt{\gamma_-}}}}{4\pi\gamma_-^{3/2}\v{x - \tilde{z}'}}
\Big(\sum_{j = 0}^{N-1}G_{j}(x-\tilde{z}')
&\Big(\frac{\sqrt{\gamma_-}}{\tau\v{x - \tilde{z}'}}\Big)^j
\label{asymptotics of the gradient of the refracted wave}
+\tilde{G}_N(x, z'; \tau)\Big),
\end{align}
where each $ G_{j}$ is a $C^\infty$ function in 
$\Bbb R^3_-$, and $\tilde{G}_N(x, z'; \tau)$ is continuous for $x \in \Bbb R^3_-$ and $z' \in \Bbb R^2$, 
and satisfies
$$
\v{\tilde{G}_N(x, z'; \tau)} \leq C_N\Big(\frac{\sqrt{\gamma_-}}{\tau\v{x - \tilde{z}'}}\Big)^N 
\qquad(x \in \Bbb R^3_-, z' \in \Bbb R^2, \tau > 1).
$$
The main task of Section \ref{Asymptotics of the refracted waves} is to show 
(\ref{asymptotics of the refracted wave}) and (\ref{asymptotics of the gradient of the refracted wave}). 

\par

From (\ref{asymptotics of the refracted wave}), 
(\ref{asymptotics of the gradient of the refracted wave}) and (\ref{expression of the FS for x_3 < 0}), 
the original problem can be reduced to finding asymptotics of the Laplace integral of the form:
\begin{align}
I(\tau; x, y) = \int_{\Bbb R^2}e^{-{\tau}l_{x, y}(z')}a(z')dz',
\qquad(x \in \overline{D}, y \in \overline{B})
\label{the Laplace integral for the refracted part}
\end{align}
where $ a \in {\mathcal B}^\infty(\Bbb R^2)$, i.e., the function $a $ belongs to 
the space of all $C^\infty$ functions
in $\Bbb R^2$ of whose all derivatives $\partial_{z'}^{\alpha'}a$ are 
bounded functions in $\Bbb R^2$. 
For $N \in {\bf N}{\cup}\{0\}$, we put 
$\V{a}_{N, {\mathcal B}^{\infty}(\Bbb R^2)} = \max_{\v{\beta} \leq N}\sup_{z' \in \Bbb R^2}
\v{\partial_{z'}^{\beta}a(z')}$. 
\par
By usual Laplace's method, the main part of the asymptotics for 
(\ref{the Laplace integral for the refracted part}) is given by
points $z'(x, y) \in \Bbb R^2$ attaining the minimum $l(x, y)$
of $l_{x, y}(z')$. We can check the point $z'(x, y)$ uniquely exists, which
corresponds to Snell's law in geometrical optics. 
Further, we can show that
$z'(x, y)$ is a $C^\infty$ function for $(x, y) \in \Bbb R^3_-\times\Bbb R^3_+$
and ${\rm Hess}(l_{x, y})(z'(x, y))$ is positive definite, where
${\rm Hess}(l_{x, y})(z') = (\partial_{z_i}\partial_{z_j}l_{x, y}(z'))$
(cf. Lemma \ref{Snell's law}). 
We put $H(x, y) = {\rm Hess}(l_{x, y})(z'(x, y))$ and
\begin{align*}
\Psi(z') = l_{x, y}(z') - l_{x, y}(z'(x, y))
- \big((H(x, y))^{-1}(z' - z'(x, y)), (z' - z'(x, y))\big).
\end{align*}
\par

Take $\phi \in C^\infty_0(\Bbb R^2)$ with $0 \leq \phi \leq 1$ and
$\phi = 1$ near the set $\{ z'(x, y) \in \Bbb R^2\vert x \in \overline{D}, y \in \overline{B}\}$, 
and divide (\ref{the Laplace integral for the refracted part}) into two parts,
\begin{align}
I(\tau; x, y) = \int_{\Bbb R^2}e^{-{\tau}l_{x, y}(z')}\phi(z')a(z')dz'
+ \int_{\Bbb R^2}e^{-{\tau}l_{x, y}(z')}(1-\phi(z'))a(z')dz'.
\label{decomposition of I(tau; x, y)}
\end{align}
Note that usual Laplace's method (cf. Theorem 7.7.5 of H\"ormander \cite{H}, which is for oscillatory integrals,
however, the proof also works for this case)  
can be applied for the first integral of (\ref{decomposition of I(tau; x, y)}). 
For the second integral of (\ref{decomposition of I(tau; x, y)}), 
integration by parts implies that this term is negligible.
Hence we obtain
\begin{align}
I(\tau; x, y) = \frac{2\pi{e}^{-{\tau}l(x, y)}}{\tau\sqrt{{\rm det }H(x, y)}}
\Big(\sum_{j = 0}^{N}(L_ja)(z'(x, y)){\tau^{-j}}
+ R_{N+1}(x, y, \tau)\Big)
\quad
(\tau > 1),
\label{asymptotic formula of I(tau; x, y)}
\end{align}
where $L_j$ is a differential operator of order less than or equal to $6j$ given by 
\begin{align*}
(L_ja)(z') = \sum_{l-k = j, 2l \geq 3k}&\frac{1}{l!k!}
\frac{(-1)^k}{2^l}
((H(x, y))^{-1}\partial_{z'}, \partial_{z'})^l((\Psi(z'))^ka))(z'),
\end{align*}
and for any $N \in {\bf N}\cup\{0\}$, there exists a constant $C_N > 0$ 
depending also on $a$ such that
$$
\v{R_{N+1}(x, y, \tau)} \leq C_N\V{a}_{2N, {\mathcal B}^{\infty}(\Bbb R^2)}\tau^{-(N+1)} 
\quad (x \in \overline{D}, y \in \overline{B},
\tau \geq 1).
$$
\par
From (\ref{asymptotic formula of I(tau; x, y)}), we can obtain 
the asymptotic expansion of $\nabla_x\Phi_\tau(x, y)$ of the form:

\begin{Prop}\label{Asymptotics of the refracted part of the gradient of the FS}
Assume $\gamma_+ < \gamma_-$. Then $\nabla_x^k\Phi_\tau(x, y)$ $(k = 0, 1)$ have
the following asymptotics:
\begin{align*}
\nabla_x^k\Phi_\tau(x, y) = 
\frac{{e}^{-{\tau}l(x, y)}}{8\pi\gamma_+\gamma_-
\sqrt{{\rm det}H(x, y)}}
\Big(\frac{-\tau}{\sqrt{\gamma_-}}\Big)^k
\Big(\sum_{j = 0}^{N}{\tau^{-j}}\Phi_{j}^{(k)}(x, y) 
+ Q_{N, \tau}^{(k)}(x, y)\Big),
\end{align*}
where $\Phi_{j}^{(k)}(x, y)$ $(k = 0, 1)$ are 
$C^\infty$ in $\overline{D}\times\overline{B}$,
for any $N \in {\rm N}\cup\{0\}$, $Q_{N, \tau}^{(k)}(x, y)$ $(k = 0, 1)$
are continuous in $\overline{D}\times\overline{B}$
with a constant $C_N > 0$ satisfying
$$
\v{Q_{N, \tau}^{(0)}(x, y)} + \v{Q_{N, \tau}^{(1)}(x, y)} \leq C_N\tau^{-(N+1)} 
\qquad(x \in \overline{D}, y \in \overline{B}, \tau \geq 1).
$$
Moreover, $\Phi_{0}^{(k)}(x, y) $ $(k = 0, 1)$ are given by
\begin{align}
\Phi_{0}^{(0)}(x, y)
&= \frac{E_{0}(x - z'(x, y))}{\v{x - \tilde{z}'(x, y)}\v{\tilde{z}'(x, y)-y}},
\label{the form of Phi_{0}^{(0)}} 
\intertext{and}
\Phi_{0}^{(1)}(x, y) 
&= \Phi_{0}^{(0)}(x, y)\frac{x - \tilde{z}'(x, y)}{\v{x - \tilde{z}'(x, y)}}.
\label{the form of Phi_{1}^{(0)}} 
\end{align}
\end{Prop}
Note that in (\ref{the form of Phi_{0}^{(0)}}), $E_0$ is the function appeared in 
(\ref{asymptotics of the refracted wave}).
The form of $E_0$ is given by (\ref{the form of E_0}) in Section \ref{Asymptotics of the refracted waves}.

\par

Proposition \ref{Asymptotics of the refracted part of the gradient of the FS}
is crucial to obtain Theorem \ref{Theorem 1.2}. 
A proof of Proposition \ref{Asymptotics of the refracted part of the gradient of the FS}
is given in Section \ref{Snell's law and the asymptotics of Phi_tau and nabla_xPhi_tau}. 
In Section \ref{Proof of Theorem ref{Theorem 1.2}}, we show Theorem \ref{Theorem 1.2}
by using Proposition \ref{Asymptotics of the refracted part of the gradient of the FS}.
This is the outline of this paper.

\setcounter{equation}{0}
\section{The refracted part of the fundamental solution}
\label{The refracted part of the fundamental solution}

In what follows, we only treat the case $y = (y', y_3)$, $y' = (y_1, y_2)$, $y_3 > 0$ for large $\tau > 0$. 
Note that the fundamental solution $\Phi_\tau(x, y)$ is given by 
$$
\Phi_\tau(x, y) = 
\begin{cases}
v_0(x, y; \tau)+v_+(x, y; \tau) & (x_3 > 0), \\
v_-(x, y; \tau) & (x_3 < 0),
\end{cases}
$$
where $v_0 = v_0(x, y; \tau)$ is the solution of 
\begin{align*}
\nabla_x\cdot(\gamma_+\nabla_x v_0) -\tau^2v_0+\delta(x-y) = 0 
\qquad \text{in } \Bbb R^3,
\end{align*}
and $v_\pm = v_\pm(x, y; \tau)$ satisfy
\begin{align}
\left\{\hskip-4mm
\begin{array}{lll}
&\nabla_x\cdot(\gamma_\pm\nabla_xv_\pm) -\tau^2v_\pm = 0
&\qquad(\pm{x_3} > 0), \\
&v_0+v_+ = v_-, \quad 
\gamma_+\partial_{x_3}v_0+\gamma_+\partial_{x_3}v_+ = \gamma_-\partial_{x_3}v_-
&\qquad(x_3 = 0).
\end{array}
\right.
\label{equation of compensation}
\end{align}
In what follows, we also write 
$x = (x', x_3)$, $x' = (x_1, x_2)$.
\par

For $v_0$, Fourier transform implies
$$
v_0(x, y; \tau) = \frac{1}{(2\pi)^3}
\int_{\Bbb R^2}d\xi'e^{i\xi'\cdot(x'-y')}
\int_{\Bbb R}
\frac{e^{i\eta\cdot(x_3-y_3)}}{\gamma_+(\xi')^2+\tau^2+\gamma_+\eta^2}
d{\eta}.
$$
We put 
$$
C_\pm^\tau(\xi') = \left(\frac{\tau^2}{\gamma_\pm}+(\xi')^2\right)^{1/2} > 0.
$$
Since 
\begin{equation}
\frac{1}{2\pi}
\int_{\Bbb R}
\frac{e^{i\eta\cdot(x_3-y_3)}}{\gamma_+(\xi')^2+\tau^2+\gamma_+\eta^2}
d{\eta}
=
\frac{1}{2\gamma_+C_+^\tau(\xi')}
e^{-C_+^\tau(\xi')\v{x_3-y_3}}, 
\label{basic integral-1}
\end{equation}
we obtain
\begin{equation*}
v_0(x, y; \tau) = \frac{1}{(2\pi)^3}
\int_{\Bbb R^2}d\xi'e^{i\xi'\cdot(x'-y')}
\frac{1}{2\gamma_+C_+^\tau(\xi')}
e^{-C_+^\tau(\xi')\v{x_3-y_3}},
\end{equation*}
which is the representation by the partial Fourier transform
\begin{equation}
\hat{v}_0(x_3, y; \xi', \tau) = 
\int_{\Bbb R^2}
e^{-i\xi'\cdot{x'}}v_0(x', x_3, y; \tau)dx'
= \frac{e^{-i\xi'{\cdot}y'}}{2\gamma_+C_+^\tau(\xi')}
e^{-C_+^\tau(\xi')\v{x_3-y_3}}
\label{the partial FT of FS for the free case}
\end{equation}
of $v_0(x, y ; \tau)$ for the tangential direction $x' \in \Bbb R^2$.
\par
To obtain $v_\pm$, we take the partial Fourier transform 
$$
\hat{v}_\pm(x_3, y; \xi', \tau) = 
\int_{\Bbb R^2}
e^{-i\xi'\cdot{x'}}v_\pm(x', x_3, y; \tau)dx',
$$
which satisfy the partial Fourier transform of (\ref{equation of compensation}), that is, 
\begin{align}
\left\{\hskip-4mm
\begin{array}{lll}
&-(\gamma_\pm(\xi')^2+\tau^2)\hat{v}_\pm+\gamma_\pm\partial_{x_3}^2\hat{v}_\pm = 0
&\qquad({\pm}x_3 > 0), \\
&\hat{v}_0+\hat{v}_+ = \hat{v}_-, \quad 
\gamma_+(\partial_{x_3}\hat{v}_0+
\partial_{x_3}\hat{v}_+) = 
\gamma_-\partial_{x_3}\hat{v}_-
&\qquad(x_3 = 0).
\end{array}
\right.
\label{equation of FT of compensation}
\end{align}
Since $v_\pm$ are bounded, from (\ref{the partial FT of FS for the free case}), 
the solutions of (\ref{equation of FT of compensation}) are given by
\begin{align}
\hat{v}_\pm(x_3, y; \xi', \tau) 
&= \frac{e^{{\mp}C^\tau_{\pm}(\xi')x_3}}{2\gamma_{\pm}C_\pm^\tau(\xi')}R_\pm^\tau(\xi')
\hat{v}_0(0; \xi', y, \tau), 
\label{FT of v_pm}
\end{align}
where
$$
R_+^\tau(\xi') = 2\gamma_{+}C_+^\tau(\xi')
\frac{\gamma_+C^\tau_+(\xi')-\gamma_-C^\tau_-(\xi')}{\gamma_+C^\tau_+(\xi')+\gamma_-C^\tau_-(\xi')}, 
\qquad
R_-^\tau(\xi') = 
\frac{4\gamma_+\gamma_-C^\tau_+(\xi')C^\tau_-(\xi')}
{\gamma_+C^\tau_+(\xi')+\gamma_-C^\tau_-(\xi')}.
$$

\par
We concentrate on for $v_-$. 
From (\ref{FT of v_pm}), (\ref{basic integral-1}) and (\ref{the partial FT of FS for the free case}), 
it follows that
$$
\hat{v}_-(x_3, y; \xi', \tau) 
= \frac{1}{2\pi}\int_{\Bbb R}
\frac{e^{i\xi_3\cdot{x_3}}}{\gamma_-\xi^2+\tau^2}d\xi_3R_-^\tau(\xi')
\hat{v}_0(0; \xi', y, \tau)
\qquad(x_3 < 0),
$$
which yields
\begin{align*}
{v}_-(x, y ;\tau) = \frac{1}{(2\pi)^6}\int_{\Bbb R^3}d\xi\int_{\Bbb R^2}dz'\int_{\Bbb R^3}d\zeta
\frac{e^{i\{\xi'\cdot(x' - z')+\xi_3\cdot{x_3}\}}}{\gamma_-\xi^2+\tau^2}R_-^\tau(\xi')
\frac{e^{i\{\zeta'\cdot(z' - y')-\zeta_3\cdot{y_3}\}}}{\gamma_+\zeta^2+\tau^2}. 
\end{align*}
Since $\Phi_\tau(x, y) = v_-(x, y; \tau)$ for $x_3 < 0$, and
$R_-^\tau(\tau\xi') = {\tau}R_-^1(\xi')$, we obtain
\begin{align*}
\Phi_\tau(x, y)
&= \frac{\tau^3}{(2\pi)^6}\int_{\Bbb R^3}{d\xi}\int_{\Bbb R^2}{dz'}
\int_{\Bbb R^3}d\zeta
\frac{e^{i\tau\{\xi'\cdot{(x'-z')}+\xi_3\cdot{x_3}\}}
}{\gamma_-\xi^2+1}
R_-^1(\xi')
\frac{e^{i\tau\{\zeta'\cdot(z'-y')-\zeta_3\cdot{y_3}\}}
}{\gamma_+\zeta^2+1}
\end{align*}
for $x_3 < 0$. Since (\ref{the form of R}) implies $R_-^1(\xi') = R_-(\v{\xi'})$,
noting (\ref{Fourier expression of FS for the free case}), (\ref{the refracted wave}) 
and the above 
formula of $\Phi_\tau(x, y)$ for $x_3 < 0$, 
we obtain (\ref{expression of the FS for x_3 < 0}).

\par

We can also obtain the formula of $\Phi_\tau(x, y)$ for $x_3 > 0$, which
is for the reflected phenomena. 
In this case, 
$\Phi_\tau(x, y) = v_0(x, y; \tau)+v_-(x, y; \tau)$, which yields
\begin{align*}
\Phi_\tau(x, y)
&= 
\frac{\tau}{(2\pi)^3}\int_{\Bbb R^3}{d\xi}
\frac{e^{i\tau\xi\cdot{(x-y)}}}{\gamma_+\xi^2+1}
\\
&
\,\,\,
+\frac{\tau^3}{(2\pi)^6}\int_{\Bbb R^3}{d\xi}\int_{\Bbb R^2}{dz'}
\int_{\Bbb R^3}d\zeta
\frac{e^{i\tau\{\xi'\cdot{(x'-z')}+\xi_3\cdot{x_3}\}}}{\gamma_+\xi^2+1}
R_+^1(\xi')
\frac{e^{i\tau\{\zeta'\cdot(z'-y')-\zeta_3\cdot{y_3}\}}}{\gamma_+\zeta^2+1},
\end{align*}
for $x_3 > 0$ similarly. In this paper, we do not use this formula.

\setcounter{equation}{0}
\section{Asymptotics of the refracted waves}
\label{Asymptotics of the refracted waves}

In this section, we show the asymptotics 
(\ref{asymptotics of the refracted wave}) and 
(\ref{asymptotics of the gradient of the refracted wave}) for the refracted wave
defined by (\ref{the refracted wave}).
We put $R(\rho) = R_-(\rho/\sqrt{\gamma_-})$ and $\tilde\tau = \tau/\sqrt{\gamma_-}$. 
From (\ref{the refracted wave}),
it follows that
$$
E^{\gamma_-}_{\tau}(x, z')
= \frac{\tau}{(2\pi)^3\gamma_-^{3/2}}\lim_{\varepsilon \to +0}
J_\varepsilon(x-\tilde{z}')
\qquad(x_3 < 0),
$$
where
$$
J_\varepsilon(x-\tilde{z}') = \int_{\Bbb R^3}e^{-\varepsilon\v{\xi'}^2}
{e}^{i\tilde\tau\xi\cdot(x-\tilde{z}')}
\frac{1}{\xi^2+1}R(\v{\xi'})d{\xi}.
$$
Note that from (\ref{basic integral-1}), it follows that
\begin{align*}
J_\varepsilon(x-\tilde{z}') &= \int_{\Bbb R^2}e^{-\varepsilon\v{\xi'}^2}
{e}^{i\tilde\tau\xi'\cdot(x'-z')}\int_{\Bbb R}{e}^{i\tilde\tau\xi_3\cdot{x_3}}
\frac{1}{\xi^2+1}d{\xi_3}R(\v{\xi'})d\xi'
\\&
= {\pi}\int_{\Bbb R^2}e^{-\varepsilon\v{\xi'}^2}
{e}^{i\tilde\tau\xi'\cdot(x'-z')}R(\v{\xi'})
e^{-\tilde\tau\v{x_3}\sqrt{1+\v{\xi'}^2}}\frac{d\xi'}{\sqrt{1+\v{\xi'}^2}},
\end{align*}
which yields
\begin{align*}
J_\varepsilon(x-\tilde{z}')
= {\pi}\int_{\Bbb R^2}e^{-\varepsilon\v{\zeta'}^2}
{e}^{i\tilde\tau\v{x'-z'}\zeta_1}R(\v{\zeta'})
e^{-\tilde\tau\v{x_3}\sqrt{1+\v{\zeta'}^2}}\frac{d\zeta'}{\sqrt{1+\v{\zeta'}^2}}
\end{align*}
by rotating the coordinate. Thus, we obtain
\begin{equation*}
E^{\gamma_-}_{\tau}(x, z')
= \frac{\tau}{2(2\pi)^2\gamma_-^{3/2}}
\int_{\Bbb R}d\zeta_2\int_{\Bbb R}{e}^{i\tilde\tau\v{x'-z'}\zeta_1}
e^{-\tilde\tau\v{x_3}\sqrt{1+\v{\zeta'}^2}}
R(\v{\zeta'})\frac{d{\zeta_1}}{\sqrt{1+\v{\zeta'}^2}}.
\end{equation*}

\par

We change the variable $\zeta_1 = \sqrt{1+\zeta_2^2}\tilde\zeta_1$, and have
\begin{align}
E^{\gamma_-}_{\tau}(x, z')
= \frac{\tau}{2(2\pi)^2\gamma_-^{3/2}}
\int_{\Bbb R}d\zeta_2\int_{\Bbb R}&{e}^{i\tilde\tau\v{x'-z'}\sqrt{1+\zeta_2^2}\zeta_1}
e^{-\tilde\tau\v{x_3}\sqrt{1+\zeta_1^2}\sqrt{1+\zeta_2^2}}
\label{For asymptotics of the refracted part 1}
\\&
{\times}
R\left(\sqrt{\zeta_1^2+\zeta_2^2+\zeta_1^2\zeta_2^2}\right)\frac{d{\zeta_1}}{\sqrt{1+\zeta_1^2}},
\nonumber
\end{align}
which yields
\begin{align}
\nabla_xE^{\gamma_-}_{\tau}(x, z')
= \frac{\tau^2}{2(2\pi)^2\gamma_-^{2}}&
\int_{\Bbb R}d\zeta_2\int_{\Bbb R}{e}^{i\tilde\tau\v{x'-z'}\sqrt{1+\zeta_2^2}\zeta_1}
e^{-\tilde\tau\v{x_3}\sqrt{1+\zeta_2^2}\sqrt{1+\zeta_1^2}}
\label{For asymptotics of the refracted part 2}
\\&\hskip-5mm
{\times}
\sqrt{1+\zeta_2^2}
R\left(\sqrt{\zeta_1^2+\zeta_2^2+\zeta_1^2\zeta_2^2}\right)
\begin{pmatrix}i\zeta_1\displaystyle\frac{x'-z'}{\v{x'-z'}} \\[5mm] 
-\frac{x_3}{\v{x_3}}\sqrt{1+\zeta_1^2}
\end{pmatrix}
\frac{d{\zeta_1}}{\sqrt{1+\zeta_1^2}}.
\nonumber
\end{align}

\par
For $x \in \Bbb R^3$, $z' \in \Bbb R^2$, we put 
\begin{align}
I_{\tilde\tau, k}(x-\tilde{z}', \zeta_2) = \int_{\Bbb R}{e}^{-\tilde\tau\sqrt{1+\zeta_2^2}
(-i\v{x'-z'}\zeta_1+\v{x_3}\sqrt{1+\zeta_1^2})}Q_k(\zeta_1, \zeta_2)
\frac{d\zeta_1}{\sqrt{1+\zeta_1^2}},
\label{the integral in zeta_1}
\end{align}
where 
\begin{align}
\left\{
\begin{array}{lll}
Q_0(\zeta_1, \zeta_2) = R\left(\sqrt{\zeta_1^2+\zeta_2^2+\zeta_1^2\zeta_2^2}\right), 
\quad
\tilde{Q}_0(\zeta_1, \zeta_2) = \sqrt{1+\zeta_2^2}Q_0(\zeta_1, \zeta_2), 
\\
Q_1(\zeta_1, \zeta_2) = i\zeta_1\tilde{Q}_0(\zeta_1, \zeta_2), 
\quad
Q_2(\zeta_1, \zeta_2) = -\sqrt{1+\zeta_1^2}\tilde{Q}_0(\zeta_1, \zeta_2). 
\end{array}
\right.
\label{definitions of Q_0 et. al.}
\end{align}

%

\par
To obtain the asymptotics of $E^{\gamma_-}_{\tau}(x, z')$ and
$\nabla_xE^{\gamma_-}_{\tau}(x, z')$, we need to study the asymptotics of 
(\ref{the integral in zeta_1}). We use the steepest decent method, which is similar to
getting the distribution kernel for the usual wave equations in the two dimensional
half-space by Hankel functions (cf. \cite{Achenbach}, p. 286 for example). 

\par
We take $\theta$ satisfying 
\begin{equation}
\sin\theta = \frac{\v{x' - z'}}{\v{x - \tilde{z}'}}, 
\quad
\cos\theta = \frac{\v{x_3}}{\v{x - \tilde{z}'}}
\quad(0 \leq \theta \leq \pi/2),
\label{the definition of theta}
\end{equation}
and put $r = \tilde\tau\v{x - \tilde{z}'}\sqrt{1+\zeta_2^2}$ and
\begin{align}
\lambda = 
\lambda(\zeta_1, x, z') = -i\sin\theta\zeta_1+\cos\theta\sqrt{1+\zeta_1^2}.
\label{change of variables}
\end{align}
Then, (\ref{the integral in zeta_1}) is written by 
\begin{align}
I_{\tilde\tau, k}(x-\tilde{z}', \zeta_2) = \int_{\Bbb R}{e}^{-r\lambda}
Q_k(\zeta_1, \zeta_2)
\frac{d\zeta_1}{\sqrt{1+\zeta_1^2}}.
\label{the integral in zeta_1'}
\end{align}
From (\ref{change of variables}), it follows that
$\zeta_1 = i\lambda\sin\theta\pm\sqrt{\lambda^2-1}\cos\theta$. Hence, putting 
$\lambda = \sqrt{1+\rho^2}$ for $\lambda \geq 1$, we have 
$\sqrt{\lambda^2-1} = \sqrt{\rho^2} = \v{\rho}$, which yields
\begin{align}
\zeta_1 = \zeta_1(\rho, x, z') = i\sqrt{1+\rho^2}\sin\theta+\rho\cos\theta
\qquad(\rho \in \Bbb R, x \in \Bbb R^3_-, z' \in \Bbb R^2).
\label{parametrization of Im lambda = 0}
\end{align}
We denote by $\Gamma $ the curve defined by (\ref{parametrization of Im lambda = 0}).  
This is the steepest decent curve of integral (\ref{the integral in zeta_1'}). 
The contour of (\ref{the integral in zeta_1'}) should be changed for $\Gamma$. 
\par
We take any $\varepsilon_0$ with $0 < \varepsilon_0 < \pi/2$. Then, for $\zeta_1 \in {\rm \bf C}$ with
$\v{\arg\zeta_1} < \pi/2-\varepsilon_0$ and $\v{\zeta_1} \geq (\sin\varepsilon_0)^{-1/2}$, 
$$
\v{\arg(1+\zeta_1^2)} \leq 2\v{\arg\zeta_1}+\v{\arg(1+\zeta_1^{-2})}
\leq \pi-2\varepsilon_0+\varepsilon_0 = \pi-\varepsilon_0, 
$$
since 
$\v{\arg(1+\zeta_1^{-2})} \leq \varepsilon_0$ for $\v{\zeta_1} \geq (\sin\varepsilon_0)^{-1/2}$.
Hence, we have
\begin{align}
\sqrt{1+\zeta_1^2} = \v{1+\zeta_1^2}^{1/2}&e^{i{\arg\zeta_1}}
e^{i\arg(1+\zeta_1^{-2})/2}
= \zeta_1(1+O(\zeta_1^{-2})) 
\label{asymptotics of sqrt{1+zeta_1^2} for Re zeta_1 > 0}
\\&\qquad
(\v{\zeta_1} \to \infty 
\text{ uniformly for } \v{\arg\zeta_1} \leq \pi/2-\varepsilon_0)
\nonumber
\end{align}
since $\sqrt{1+\zeta_1^2} > 0$ for $\zeta_1 \in \Bbb R$. 
Similarly, we also obtain
\begin{align}
\sqrt{1+\zeta_1^2} = \v{1+\zeta_1^2}^{1/2}&e^{i(2{\arg\zeta_1}+2\pi)/2}
e^{i\arg(1+\zeta_1^{-2})/2}
= e^{\pi{i}}\zeta_1(1+O(\zeta_1^{-2})) 
\label{asymptotics of sqrt{1+zeta_1^2} for Re zeta_1 < 0}\\&\qquad
(\v{\zeta_1} \to \infty 
\text{ uniformly for } \v{\arg\zeta_1+\pi} \leq \pi/2-\varepsilon_0)
\nonumber
\end{align}
since in this case, it follows that 
$$
\v{\arg(1+\zeta_1^2)+2\pi} \leq 2\v{\arg\zeta_1+\pi}
+\v{\arg(1+\zeta_1^{-2})}
\leq \pi-2\varepsilon_0 + \varepsilon_0 = \pi-\varepsilon_0.
$$
From these asymptotics, it follows that $\lambda$ defined by 
(\ref{change of variables}) satisfies 
\begin{align}
{\rm Re }\lambda &= {\rm Im }\zeta_1\sin\theta+{\rm Re }\zeta_1\cos\theta
+O(\v{\zeta_1}^{-1})
\label{asymptotics of lambda in right half space}
\\&\hskip30mm
(\v{\zeta_1} \to \infty 
\text{ uniformly for } \v{\arg\zeta_1} \leq \pi/2-\varepsilon_0),
\nonumber
\\
{\rm Re }\lambda &= {\rm Im }\zeta_1\sin\theta-{\rm Re }\zeta_1\cos\theta
+O(\v{\zeta_1}^{-1})
\label{asymptotics of lambda in left half space}\\&\hskip30mm
(\v{\zeta_1} \to \infty 
\text{ uniformly for } \v{\arg\zeta_1+\pi} \leq \pi/2-\varepsilon_0).
\nonumber
\end{align}
Noting $1+\zeta_1^2 = (\sqrt{1+\rho^2}\cos\theta+i\rho\sin\theta)^2$, and 
(\ref{asymptotics of sqrt{1+zeta_1^2} for Re zeta_1 > 0}) and
(\ref{asymptotics of sqrt{1+zeta_1^2} for Re zeta_1 < 0}), we also have
\begin{align}
\sqrt{1+\zeta_1^2} = \sqrt{1+\rho^2}\cos\theta+{i}\rho\sin\theta
\label{expression of sqrt{1+zeta_1^2} on Im lambda = 0}
\end{align}
for $\zeta_1 = \zeta_1(\rho, x, z')$.

\par

From (\ref{the form of R}), it follows that
$$
R(\rho) = R_-(\rho/\sqrt{\gamma_-})
= \frac{4\sqrt{\gamma_-}\sqrt{a_0^2 + \rho^2}\sqrt{1 + \rho^2}}
{\sqrt{a_0^2 + \rho^2}+a_0^2\sqrt{1 + \rho^2}},
$$
where 
\begin{equation}
\displaystyle
a_0 = \sqrt{\frac{\gamma_-}{\gamma_+}}.
\label{definition of a_0}
\end{equation} 
Since $1+(\sqrt{\zeta_1^2+\zeta_2^2+\zeta_1^2\zeta_2^2})^2 = (1+\zeta_1^2)(1+\zeta_2^2)$, 
we have
\begin{align}
R\left(\sqrt{\zeta_1^2+\zeta_2^2+\zeta_1^2\zeta_2^2}\right)
= \frac{4\sqrt{\gamma_-}\sqrt{1+\zeta_2^2}\sqrt{1 + \zeta_1^2}P(\zeta_1, \zeta_2)}
{P(\zeta_1, \zeta_2)+a_0^2\sqrt{1 + \zeta_1^2}}, 
\label{the form of composition function of R}
\end{align}
where
$$
P(\zeta_1, \zeta_2) = \sqrt{\frac{a_0^2-1}{1+\zeta_2^2}+1+\zeta_1^2}. 
$$

\par

In what follows, we assume $\gamma_+ < \gamma_-$, being the case that 
there is no total reflected wave for incident waves coming from the lower side 
(cf. Remark \ref{Why gamma_+ < gamma_-}). 
In this case, $a_0 > 1$. Hence, 
$P(\zeta_1, \zeta_2)$ and $\sqrt{1+\zeta_1^2}$ are holomorphic for 
$\zeta_1 \in {\rm {\bf C}}\setminus((-i\infty, -i]\cup[i, i\infty))$. From this, 
(\ref{asymptotics of lambda in right half space}) and (\ref{asymptotics of lambda in left half space}),
we can change the contour of (\ref{the integral in zeta_1'}) for $\Gamma$, which yields
\begin{align*}
I_{\tilde\tau, k}(x-\tilde{z}', \zeta_2) = \int_{\Gamma}{e}^{-r\lambda}
Q_k(\zeta_1, \zeta_2)
\frac{d\zeta_1}{\sqrt{1+\zeta_1^2}}.
\end{align*}
We can express this integral by using 
the parametrization of $\Gamma$ given by (\ref{parametrization of Im lambda = 0}).
In this parametrization, $\lambda = \sqrt{1+\rho^2}$, and 
(\ref{parametrization of Im lambda = 0}) and 
(\ref{expression of sqrt{1+zeta_1^2} on Im lambda = 0}) implies
$$\displaystyle
\frac{d\zeta_1}{d\rho} = \frac{\sqrt{1+\zeta_1^2}}{\sqrt{1+\rho^2}},
$$
and we obtain 
\begin{align}
I_{\tilde\tau, k}(x-\tilde{z}', \zeta_2) = 
\int_{\Bbb R}{e}^{-{\tilde\tau}\v{x - \tilde{z}'}\sqrt{1+\zeta_2^2}\sqrt{1+\rho^2}}
Q_k(\zeta_1(\rho, x, z'), \zeta_2)\frac{d\rho}{\sqrt{1+\rho^2}}.
\label{parametrized expression for I_k}
\end{align}
Note that $Q_k(\zeta_1(\rho, x, z'), \zeta_2)$ are $C^\infty$ function of $\rho \in \Bbb R$
since $a_0 > 1 $ imples $\Gamma \subset {\rm {\bf C}}\setminus((-i\infty, -i]\cup[i, i\infty))$. 
For simplicity, we write $\sigma_1 = \rho$, $\sigma_2 = \zeta_2$ and $\sigma = (\sigma_1, \sigma_2)$,
and put 
\begin{align}
f(\sigma) &= \sqrt{1+\sigma_1^2}\sqrt{1+\sigma_2^2},
\quad
F_k(\sigma, x, z') = Q_k(\zeta_1(\sigma_1, x, z'), \sigma_2)
\frac{1}{\sqrt{1+\sigma_1^2}}
\nonumber
\\
F_+(\sigma, x, z') &= \sqrt{1+\sigma_1^2}\tilde{Q}_0(\zeta_1(\sigma_1, x, z'), \sigma_2)
\frac{1}{\sqrt{1+\sigma_1^2}}, 
\qquad
\nonumber
\\[-4mm]
\intertext{and} 
\nonumber\\[-10mm]
F_-(\sigma, x, z') &= i\sigma_1\tilde{Q}_0(\zeta_1(\sigma_1, x, z'), \sigma_2)
\frac{1}{\sqrt{1+\sigma_1^2}}. 
\label{definitions of F_-}
\end{align}
Notice that (\ref{definitions of Q_0 et. al.}), 
(\ref{parametrization of Im lambda = 0}) and 
(\ref{expression of sqrt{1+zeta_1^2} on Im lambda = 0}) imply that
\begin{align}
\left\{
\begin{array}{lll}
F_{1}(\sigma, x, z') = F_-(\sigma, x, z')\cos\theta - F_+(\sigma, x, z')\sin\theta, 
\\
F_{2}(\sigma, x, z') = -F_+(\sigma, x, z')\cos\theta - F_-(\sigma, x, z')\sin\theta.
\end{array}
\right.
\label{relations between F_1, F_2 and F_pm}
\end{align}

From (\ref{For asymptotics of the refracted part 1}), 
(\ref{For asymptotics of the refracted part 2}) 
and (\ref{parametrized expression for I_k}), it follows that
\begin{align}
E^{\gamma_-}_{\tau}(x, z')
&= \frac{\tau}{2(2\pi)^2\gamma_-^{3/2}}
\int_{\Bbb R^2}{e}^{-{\tilde\tau}\v{x - \tilde{z}'}f(\sigma)}
F_0(\sigma, x, z')d\sigma, 
\label{For asymptotics of the refracted part 1'}
\\
\nabla_{x'}E^{\gamma_-}_{\tau}(x, z') &= \frac{\tau^2}{2(2\pi)^2\gamma_-^{2}}
\int_{\Bbb R^2}{e}^{-{\tilde\tau}\v{x - \tilde{z}'}f(\sigma)}
F_1(\sigma, x, z')d\sigma\displaystyle\frac{x'-z'}{\v{x'-z'}}, 
\label{For asymptotics of the refracted part 2'-x'}
\\
\partial_{x_3}E^{\gamma_-}_{\tau}(x, z') &= \frac{\tau^2}{2(2\pi)^2\gamma_-^{2}}
\int_{\Bbb R^2}{e}^{-{\tilde\tau}\v{x - \tilde{z}'}f(\sigma)}
F_2(\sigma, x, z')d\sigma\displaystyle\frac{x_3}{\v{x_3}}.
\label{For asymptotics of the refracted part 2'-x_3}
\end{align}

\par

From (\ref{relations between F_1, F_2 and F_pm})-(\ref{For asymptotics of the refracted part 2'-x_3}), 
the problem is reduced to finding the asymptotics of 
\begin{align*}
\int_{\Bbb R^2}{e}^{-{\tilde\tau}\v{x - \tilde{z}'}f(\sigma)}
F_k(\sigma, x, z')d\sigma
\qquad(k = 0, 1, 2 \text{ and } \pm)
\end{align*}
as $\tilde\tau \to \infty$.
For treating these integrals, we need to assume $\gamma_+ < \gamma_-$, which implies that 
the amplitude functions $F_k$ in 
(\ref{For asymptotics of the refracted part 1'})-(\ref{For asymptotics of the refracted part 2'-x_3})  
are smooth. This allows to use the Laplace methods to give the asymptotic expansions for 
$E^{\gamma_-}_{\tau}(x, z')$ and its gradient.
\begin{Prop}\label{Asymptotics of refracted parts for case 1}
Assume $\gamma_+ < \gamma_-$. Then, it follows that
\begin{align*}
E^{\gamma_-}_{\tau}(x, z')
= \frac{e^{-{\tau\v{x-\tilde{z}'}/\sqrt{\gamma_-}}}}{4\pi\gamma_-\v{x - \tilde{z}'}}
&\Big(\sum_{j = 0}^{N-1}E_{j}(x-\tilde{z}')
\Big(\frac{\sqrt{\gamma_-}}{\tau\v{x - \tilde{z}'}}\Big)^j
+\tilde{E}_{N}(x, z'; \tau)
\Big),
\intertext{and for $k = 1, 2, 3$, } 
\partial_{x_k}E^{\gamma_-}_{\tau}(x, z')
= \frac{-{\tau}e^{-{\tau\v{x-\tilde{z}'}/\sqrt{\gamma_-}}}}{4\pi\gamma_-^{3/2}\v{x - \tilde{z}'}}
&\Big(
\sum_{j = 0}^{N-1}G_{k, j}(x-\tilde{z}')\Big(\frac{\sqrt{\gamma_-}}{\tau\v{x - \tilde{z}'}}\Big)^j
+\tilde{G}_{k, N}(x, z'; \tau)\Big), 
\end{align*}
where $E_j(x - \tilde{z}')$, $G_{k, j}(x - \tilde{z}')$ 
($k = 1, 2, 3$ and $j = 0, 1, 2, \ldots)$ are $C^\infty$
functions for $x \in \Bbb R^3_-$ and $z' \in \Bbb R^2$. 
Here, the remainder terms
$\tilde{E}_{N}(x, z'; \tau)$ and $\tilde{G}_{k, N}(x, z'; \tau)$ $(k = 1, 2, 3)$ 
are estimated by 
$$
\v{\tilde{E}_{N}(x, z'; \tau)}+\sum_{k = 1}^3\v{\tilde{G}_{k, N}(x, z'; \tau)}
\leq C_N\Big(\frac{\sqrt{\gamma_-}}{\tau\v{x - \tilde{z}'}}\Big)^{N}
\quad(x \in \Bbb R^3_-, z' \in \Bbb R^2)
$$
for some constant $C_N > 0$ depending only on $N \in {\rm {\bf N}}$.
In particular, we have
\begin{align}
E_{0}(x - \tilde{z}') &= \frac{4\sqrt{\gamma_-}\v{x_3}\sqrt{a_0^2\v{x - \tilde{z}'}^2- \v{x' - z'}^2}}
{\v{x - \tilde{z}'}\big(\sqrt{a_0^2\v{x - \tilde{z}'}^2- \v{x' - z'}^2}+a_0^2\v{x_3}\big)},
\label{the form of E_0}
\end{align}
where $a_0 > 1$ is given by (\ref{definition of a_0}), and
\begin{align*}
\left\{
\begin{array}{ll}
G_{k, 0}(x - \tilde{z}') = E_{0}(x - \tilde{z}')\displaystyle\frac{x_k - z_k}{\v{x - \tilde{z}'}}
\quad(k = 1, 2)\quad
\\[2mm]
G_{3, 0}(x - \tilde{z}') = E_{0}(x - \tilde{z}')\displaystyle\frac{x_3}{\v{x - \tilde{z}'}}. 
\end{array}
\right.
\end{align*}
\end{Prop}
\par\noindent
Proof: Note that $f$ in the integrals of 
(\ref{For asymptotics of the refracted part 1'})-(\ref{For asymptotics of the refracted part 2'-x_3}) 
has only one critical point $\sigma = 0$, and 
${\rm Hess}f(0) = I$, where ${\rm Hess}f(0)$ is the Hessian of $f$ 
at $\sigma = 0$ and $I$ is the $2{\times}2$ unit matrix. Since $(1+\sigma_1^2)(1+\sigma_2^2) \geq 1+\v{\sigma}^2
\geq (1+\v{\sigma}/3)^2$ for $\v{\sigma} \geq 3/4$, we have
\begin{align}
f(\sigma) \geq 1+\frac{\v{\sigma}}{3} \geq \frac{9}{8}+\frac{\v{\sigma}}{6}\quad (\v{\sigma} \geq 3/4).
\label{estimates of f and G_k for non-glancing case 1}
\end{align}
Further, there exists a constant $C > 0$ such that
\begin{align}
\v{F_k(\sigma, x, z')} \leq C(1+\v{\sigma})^3 \quad(\sigma \in \Bbb R^2, x \in \Bbb R^3_-, z' \in \Bbb R^2).
\label{estimates of f and G_k for non-glancing case 1'}
\end{align}
Take $\psi \in C^\infty_0(\Bbb R^2)$ with $0 \leq \psi \leq 1$, 
$\psi(\sigma) = 1$ for $\v{\sigma} \leq 1$,
and $\psi(\sigma) = 0$ $\v{\sigma} \geq 3/2$. 
From (\ref{estimates of f and G_k for non-glancing case 1}) and
(\ref{estimates of f and G_k for non-glancing case 1'}), it follows that
\begin{align}
\Big\vert\int_{\Bbb R^2}e^{-{\tilde\tau}\v{x - \tilde{z}'}f(\sigma)}
F_k(\sigma, x, z')(1-\psi(\sigma))d\sigma\Big\vert
&\leq Ce^{-9{\tilde\tau\v{x - \tilde{z}'}}/8}
\int_{\Bbb R^2}(1+\v{\sigma})^3e^{-({\tilde\tau}\v{x - \tilde{z}'}/6)\v{\sigma}}d\sigma
\nonumber\\
&\leq \frac{C_Ne^{-{\tilde\tau\v{x - \tilde{z}'}}}}{(\tilde\tau\v{x - \tilde{z}'})^{N}},
\label{negrigible part of the asymptotics}
\end{align}
and usual Laplace's method as is stated in Introduction implies
\begin{align}
\Big\vert e^{{\tilde\tau\v{x - \tilde{z}'}}}\int_{\Bbb R^2}e^{-{\tilde\tau}\v{x - \tilde{z}'}f(\sigma)}
F_k(\sigma,  x, z')&\psi(\sigma)d\sigma
- (\frac{2\pi}{{\tilde\tau}\v{x - \tilde{z}'}})
\sum_{j = 0}^{N-1}F_{k, j}(x - \tilde{z}')({\tilde\tau}\v{x - \tilde{z}'})^{-j}
\Big\vert
\nonumber\\&\hskip-10mm
\leq C_N\V{F_k(\cdot, x, z')}_{2N, {\mathcal B}^\infty(\overline{B_2(0)})}
(\tilde\tau\v{x - \tilde{z}'})^{-N-1},
\label{asymptotic expansion for the integral from the FS}
\end{align}
where $B_2(0) = \{\sigma \in \Bbb R^2 \vert \v{\sigma} < 2\}$ and 
$$
\V{F_k(\cdot, x, z')}_{2N, {\mathcal B}^\infty(\overline{B_2(0)})} 
= \displaystyle\max_{\v{\beta} \leq 2N}\sup_{\sigma \in B_2(0)}
\v{\partial_\sigma^{\beta}F_k(\sigma, x, z')}.
$$
In (\ref{asymptotic expansion for the integral from the FS}), 
$F_{k, j}(x - \tilde{z}')$ is given by
\begin{align}
F_{k, j}(x - \tilde{z}') &
= \sum_{l-p = j, 2l \geq 3p}\frac{1}{l!p!}
\frac{(-1)^p}{2^l}\triangle^l((\Psi^pF_k(\cdot, x, z')))(0)
\qquad(j = 0, 1, \ldots),  
\label{definition of F_{k, j}}
\end{align}
where $\Psi(\sigma) = f(\sigma) - f(0) - \v{\sigma}^2/2$.

\par
Since $(\Psi(\sigma))^p = O(\v{\sigma}^{4p})$, for 
\begin{align*}
\Psi(\sigma) = \frac{(\sigma_1\sigma_2)^2}{4}- \frac{\sigma_1^4+\sigma_2^4}{8}+O(\v{\sigma}^6) 
\qquad(\v{\sigma} \to 0),
\end{align*}
we have $\triangle^l((\Psi(\sigma))^p)(0) = 0$ for $2l < 4p$, which yields that 
$p \leq j$ holds for $ l - p = j$ with $\triangle^l((\Psi(\sigma))^p)(0) \neq 0$.
Thus, the summation of (\ref{definition of F_{k, j}}) is not taken from 
all pairs $(p, l)$ in (\ref{definition of F_{k, j}}). It consists only from 
pairs $(p, l)$ with $0 \leq p \leq j$ and $l = p+j$, 
which implies
\begin{align}
F_{k, j}(x - \tilde{z}') 
= \sum_{p = 0}^{j}\frac{1}{(p+j)!p!}\frac{(-1)^p}{2^{p+j}}
\triangle^{p+j}((\Psi^pF_k(\cdot, x, z')))(0)
\qquad(j = 0, 1, \ldots).
\label{the forms of F_{k, j}(z)}
\end{align}
From (\ref{the definition of theta}) and (\ref{relations between F_1, F_2 and F_pm}), it follows that
\begin{align}
\left\{
\begin{array}{ll}
F_{1, j}(x - \tilde{z}') &= F_{-, j}(x - \tilde{z}')\displaystyle\frac{\v{x_3}}{\v{x - \tilde{z}'}} 
- F_{+, j}(x - \tilde{z}')\frac{\v{x' - z'}}{\v{x - \tilde{z}'}}, 
\\[3mm]
F_{2, j}(x - \tilde{z}') &= -F_{+, j}(x - \tilde{z}')\displaystyle\frac{\v{x_3}}{\v{x - \tilde{z}'}}
- F_{-, j}(x - \tilde{z}')\frac{\v{x' - z'}}{\v{x - \tilde{z}'}}.
\end{array}
\right.
\label{relations between F_{pm, j} and F_{k, j}}
\end{align}

Since (\ref{the forms of F_{k, j}(z)}) and (\ref{parametrization of Im lambda = 0})
imply that for $k = 0, 1, 2$, 
$$
F_{k, 0}(x - \tilde{z}') = Q_k(\zeta_1(0, x, z'), 0) = Q_k(i\sin\theta, 0).
$$
From (\ref{definitions of Q_0 et. al.}) and (\ref{the form of composition function of R}), 
it follows that
which yields
\begin{align}
F_{0, 0}(x - \tilde{z}') &= \frac{4\sqrt{\gamma_-}\cos\theta\sqrt{a_0^2- \sin^2\theta}}
{\sqrt{a_0^2- \sin^2\theta}+a_0^2\cos\theta},
\label{form of F_{0, 0} no.1}
\\
F_{1, 0}(x - \tilde{z}') &= -F_{0, 0}(x - \tilde{z}')\sin\theta,  
\qquad
F_{2, 0}(x - \tilde{z}') = -F_{0, 0}(x - \tilde{z}')\cos\theta,
\label{forms of F_{1, 0} and F_{2, 0} no.1}
\end{align}
where the relations between $\theta$ and $x-\tilde{z}'$ is given by
(\ref{the definition of theta}). 
\par

Combining (\ref{For asymptotics of the refracted part 1'})-(\ref{For asymptotics of the refracted part 2'-x_3}) 
with (\ref{negrigible part of the asymptotics}) and
(\ref{asymptotic expansion for the integral from the FS}), 
we obtain the asymptotics stated in 
Proposition \ref{Asymptotics of refracted parts for case 1} except
the properties of the coefficients functions $E_{j}$ and $G_{k, j}$.
Notice that $E_{j}$ and $G_{k, j}$ are given by 
\begin{align}
E_j(x - \tilde{z}') = F_{0, j}(x - \tilde{z}')  
\label{formula of E_j}
\end{align}
and
\begin{align}
\left\{
\begin{array}{ll}
G_{k, j}(x - \tilde{z}') = -F_{1, j}(x - \tilde{z}')\displaystyle\frac{x_k-z_k}{\v{x'-z'}} \qquad(k = 1, 2), \\
G_{3, j}(x - \tilde{z}') = -F_{2, j}(x - \tilde{z}')\displaystyle\frac{x_3}{\v{x_3}}.
\end{array}
\right.
\label{formula of G_{k, j}}
\end{align}
From (\ref{form of F_{0, 0} no.1}), (\ref{forms of F_{1, 0} and F_{2, 0} no.1}) 
and (\ref{the definition of theta}), we also have the forms of 
$E_{0}$, $G_{1, 0}$, $G_{2, 0}$ and $G_{3, 0}$. The remaining parts are to
prove smoothness of the coefficients.

\par

For $N \in {\bf {\rm N}}\cup\{0\}$, we denote by ${\mathcal P}_N$ the set 
consisting of functions $p$ of $\theta$ of the form:
$$
p(\theta) = \sum_{j + 2k \leq 2N}a_{jk}(\sin^2\theta)
(\sin^2\theta-\cos^2\theta)^j(\cos\theta\sin\theta)^{2k},
$$
where $a_{jk}(t) $ are $C^\infty$ for $\v{t} < 1+\delta$ with some
positive $\delta > 0$. Note that any $p \in {\mathcal P}_N$ is
regarded as a $C^\infty $ function in $x - \tilde{z}' \in \Bbb R^3_-$
by relations (\ref{the definition of theta}). 

\par

First, we show smoothness of $E_{j}$. 
Since $\Psi$ is an even function with respect to each of $\sigma_1$ and $\sigma_2$, 
(i.e. 
$\Psi(-\sigma_1, \sigma_2) = \Psi(\sigma)$ and 
$\Psi(\sigma_1, -\sigma_2) = \Psi(\sigma)$), and $Q_0(\zeta_1, \zeta_2)$ is a function for $\zeta_1^2$
and $\zeta_2^2$, for $0 \leq p \leq j$, $(\Psi(\sigma))^pF_0(\sigma, x, z')$ is of the form:
\begin{align}
(\Psi(\sigma))^pF_0(\sigma, x, z') = A(\sigma, \sigma_1^2(\sin^2\theta-\cos^2\theta)+\sin^2\theta, 
\varphi(\sigma_1)\cos\theta\sin\theta),
\label{the form of Psi^pF_0}
\end{align}
where $ A = A(\sigma, \eta)$ for $\sigma, \eta \in \Bbb R^2$ is a $C^\infty$ function
for $\v{\sigma} < \delta$, $\v{\eta} < 1+\delta$ with a sufficiently small $\delta > 0$,
and $A$ is an even function with respect to each of $\sigma_1$ and $\sigma_2$, (i.e. 
$A(-\sigma_1, \sigma_2, \eta) = A(\sigma_1, -\sigma_2, \eta) = A(\sigma, \eta)$), and
$\varphi$ is a $C^\infty$ and an odd function in a neighborhood of $\sigma_1 = 0$.
From (\ref{the forms of F_{k, j}(z)}) and (\ref{formula of E_j}), 
it suffices to show for any $\v{\alpha} \leq 2j$, 
\begin{align}
\partial_\sigma^{2\alpha}\big((\Psi(\sigma))^p&F_0(\sigma, x, z')\big)\big\vert_{\sigma = 0} 
\in {\mathcal P}_{2j}. 
\label{target for F_{0, j}(z)}
\end{align}

\par

Since $A$ is an even function with respect to $\sigma_1$ and $\sigma_2$, it follows that
\begin{align}
\partial_\sigma^{2\alpha}&\big((\Psi(\sigma))^pF_0(\sigma, x, z')\big)\vert_{\sigma = 0} 
= \sum_{l = 0}^{\alpha_1}\frac{(2\alpha_1)!}{(2l)!(2(\alpha_1-l))!}
\label{delivatives of Psi^pF_0}
A_{l, \alpha_1, \alpha_2}(\theta),
\end{align}
where $A_{l, \alpha_1, \alpha_2}$ is defined by 
\begin{align*}
A_{l, \alpha_1, \alpha_2}(\theta) = 
\partial_\rho^{2l}\big\{\big(\partial_{\sigma_1}^{2\alpha_1-2l}
\partial_{\sigma_2}^{2\alpha_2}A\big)(0, \rho^2(\sin^2\theta-\cos^2\theta)+\sin^2\theta, 
\varphi(\rho)\cos\theta\sin\theta)\big\}\vert_{\rho = 0}
\end{align*}
for $0 \leq l \leq \alpha_1$ and $0 \leq \alpha_2$.   
Since $\varphi(\rho)$ is odd, 
$\partial_\rho^{2l}\big[\rho^{2\gamma_1}(\varphi(\rho))^{\gamma_2}
\big]\big\vert_{\rho = 0} = 0
$
for any odd $\gamma_2$, 
Taylor's theorem implies $A_{l, \alpha_1, \alpha_2} \in {\mathcal P}_{2l}$ since
$A_{l, \alpha_1, \alpha_2}$ can be written as
\begin{align*}
A_{l, \alpha_1, \alpha_2}(\theta) &= \sum_{\gamma_1+2\gamma_2 \leq 2l}
\frac{\big(\partial_{\eta_1}^{\gamma_1}\partial_{\eta_2}^{2\gamma_2}\partial_{\sigma_1}^{2\alpha_1-2l}
\partial_{\sigma_2}^{2\alpha_2}A\big)(0, \sin^2\theta, 0)}{\gamma_1!(2\gamma_2)!}
\partial_\rho^{2l}\big[(\rho^{2\gamma_1}(\varphi(\rho))^{2\gamma_2}
\big]\big\vert_{\rho = 0}
\\&
\hskip20mm
\times
(\sin^2\theta-\cos^2\theta)^{\gamma_1}
(\cos\theta\sin\theta)^{2\gamma_2}.
\end{align*}
Combining $A_{l, \alpha_1, \alpha_2} \in {\mathcal P}_{2l}$ shown in the above, 
(\ref{the form of Psi^pF_0}) with (\ref{delivatives of Psi^pF_0}), 
we obtain (\ref{target for F_{0, j}(z)}),
which yields 
smoothness of $E_{j}(x - \tilde{z}') $ in $x \in \Bbb R^3_-$ and $z \in \Bbb R^2$
for (\ref{the forms of F_{k, j}(z)}).

\par
Next, we show smoothness of $G_{k, j}$. From 
(\ref{relations between F_{pm, j} and F_{k, j}}) and
(\ref{formula of G_{k, j}}), 
$G_{k, j}$ are given by
\begin{align*}
G_{k, j}(x - \tilde{z}') &= -\left(F_{-, j}(x - \tilde{z}')
\displaystyle\frac{\v{x_3}}{\v{x' - z'}} 
- F_{+, j}(x - \tilde{z}')\right)\frac{x_k - z_k}{\v{x - \tilde{z}'}}
\quad(k = 1, 2),
\\
G_{3, j}(x - \tilde{z}') &= -\left(- F_{+, j}(x - \tilde{z}')
-F_{-, j}(x - \tilde{z}')
\displaystyle\frac{\v{x' - z'}}{\v{x_3}} 
\right)\frac{x_3}{\v{x - \tilde{z}'}}.
\end{align*}
Since $\tilde{Q}_0(\zeta_1, \zeta_2)$ is also a function for $\zeta_1^2$
and $\zeta_2^2$, $(\Psi(\sigma))^pF_+(\sigma, x, z')$ is the same form as (\ref{the form of Psi^pF_0}),
which yields that $F_{+, j}(x - \tilde{z}') $ is $C^\infty$ for $x \in \Bbb R^3_-$ and $z \in \Bbb R^2$.
Thus, it suffices to show that 
$F_{-, j}(x - \tilde{z}')\frac{\v{x_3}}{\v{x' - z'}} $
is $C^\infty$ for $x \in \Bbb R^3_-$ and $z \in \Bbb R^2$.

\par

From (\ref{definitions of F_-}), by using a function $A_-(\sigma, \eta)$
with the same property as for $A$
in (\ref{the form of Psi^pF_0}), we can write 
$(\Psi(\sigma))^pF_-(\sigma, x, z')$ as
\begin{align}
(\Psi(\sigma))^pF_-(\sigma, x, z') &= 
\sigma_1A_-(\sigma, \sigma_1^2(\sin^2\theta-\cos^2\theta)+\sin^2\theta, 
\varphi(\sigma_1)\cos\theta\sin\theta).
\label{the form of Psi^pF_1}
\end{align}
Since there is $\sigma_1$ in (\ref{the form of Psi^pF_1}) for $\v{\alpha} \leq j$, we have 
$$
\begin{array}{l}
\displaystyle
\,\,\,\,\,\,
\partial_\sigma^{2\alpha}
\big(
\sigma_1A_-(\sigma, \sigma_1^2(\sin^2\theta-\cos^2\theta)+\sin^2\theta, 
\varphi(\sigma_1)\cos\theta\sin\theta)\big)\vert_{\sigma = 0}
\\
\\
\displaystyle
= \sum_{l = 1}^{\alpha_1}\frac{(2\alpha_1)!}
{(2l-1)!(2\alpha_1-2l)!}
\\
\\
\displaystyle
\,\,\,
\times
\partial_\rho^{2l-1}\big\{\big(\partial_{\sigma_1}^{2\alpha_1-2l}
\partial_{\sigma_2}^{2\alpha_2}A_-\big)
(0, 0, \rho^2(\sin^2\theta-\cos^2\theta)+\sin^2\theta, 
\varphi(\rho)\cos\theta\sin\theta)\big\}\vert_{\rho = 0}.
\end{array}
$$
Moreover, Taylor's theorem implies 
$$\begin{array}{l}
\displaystyle
\,\,\,\,\,\,
\partial_\rho^{2l-1}\big\{\big(\partial_{\sigma_1}^{2\alpha_1-2l}
\partial_{\sigma_2}^{2\alpha_2}A_-\big)(0, \rho^2(\sin^2\theta-\cos^2\theta)+\sin^2\theta, 
\varphi(\rho)\cos\theta\sin\theta)\big\}\vert_{\rho = 0}
\\
\\
\displaystyle
= \sum_{\gamma_1+2\gamma_2+1 \leq 2l-1}
\frac{\big(\partial_{\eta_1}^{\gamma_1}\partial_{\eta_2}^{2\gamma_2+1}\partial_{\sigma_1}^{2\alpha_1-2l}
\partial_{\sigma_2}^{2\alpha_2}A_-\big)(0, \sin^2\theta, 0)}{\gamma_1!(2\gamma_2+1)!}
\partial_\rho^{2l-1}\big[(\rho^{2\gamma_1}(\varphi(\rho))^{2\gamma_2+1}
\big]\big\vert_{\rho = 0}
\\
\\
\displaystyle
\,\,\,
\times
(\sin^2\theta-\cos^2\theta)^{\gamma_1}
(\cos\theta\sin\theta)^{2\gamma_2+1}
\end{array}
$$
since $\partial_\rho^{2l-1}\big[(\rho^{2\gamma_1}(\varphi(\rho))^{\gamma_2}
\big]\big\vert_{\rho = 0} = 0$ for even $\gamma_2$.
From these equalities, we obtain 
\begin{align*}
\partial_\sigma^{2\alpha}&\big(
\sigma_1A_-(\sigma, \sigma_1^2(\sin^2\theta-\cos^2\theta)+\sin^2\theta, 
\varphi(\sigma_1)\cos\theta\sin\theta)\big)\vert_{\sigma = 0} 
= Y(\theta)\cos\theta\sin\theta
\end{align*}
with some $Y \in {\mathcal P}_{2j}$. From this property and (\ref{the forms of F_{k, j}(z)}),
$F_{-, j}(x - \tilde{z}')$ is of the form:
$$
F_{-, j}(x - \tilde{z}') = \tilde{F}_{-, j}(x - \tilde{z}')
\frac{\v{x_3}\v{x' - z'}}{\v{x - \tilde{z}'}^2}
$$
with some $C^\infty$ function $\tilde{F}_{-, j}(x - \tilde{z}')$ for 
$x \in \Bbb R^3_-$ and $z \in \Bbb R^2$. 
Hence, we have
$$\left\{
\begin{array}{l}
\displaystyle
F_{-, j}(x - \tilde{z}')\frac{\v{x_3}}{\v{x' - z'}}
= \tilde{F}_{-, j}(x - \tilde{z}')\frac{\v{x_3}^2}{\v{x - \tilde{z}'}^2}, \\
\\
\displaystyle
F_{-, j}(x - \tilde{z}')\frac{\v{x' - z'}}{\v{x_3}}
= \tilde{F}_{-, j}(x - \tilde{z}')\frac{\v{x' - z'}^2}{\v{x - \tilde{z}'}^2}, 
\end{array}
\right.
$$
which complete the proof of 
Proposition \ref{Asymptotics of refracted parts for case 1}.
\hfill$\blacksquare$
\vskip1pc\noindent
\par\noindent
\begin{Remark}\label{Asymptotics of gradient of refracted parts for case 1}
We put $G_j(x - \tilde{z}') = {}^t(G_{1, j}(x - \tilde{z}'), 
G_{2, j}(x - \tilde{z}'), G_{3, j}(x - \tilde{z}'))$. 
From the proof of Proposition \ref{Asymptotics of refracted parts for case 1}, 
the vector valued functions $G_j(x - \tilde{z}')$ are of the form
$$
G_{j}(x-\tilde{z}') = F_{+, j}(x - \tilde{z}')\frac{x - \tilde{z}'}{\v{x - \tilde{z}'}}
- \tilde{F}_{-, j}(x-\tilde{z}')
\begin{pmatrix}
\displaystyle
\frac{\v{x_3}^2}{\v{x - \tilde{z}'}^2}\frac{x'- z'}{\v{x-\tilde{z}'}} 
\\[3mm]
\displaystyle
-\frac{\v{x' - z'}^2}{\v{x - \tilde{z}'}^2}\frac{x_3}{\v{x-\tilde{z}'}}
\end{pmatrix}.
$$
Using these $G_j$ and putting 
$\tilde{G}_N(x, z'; \tau)
= {}^t(\tilde{G}_{1, N}(x, z'; \tau), \tilde{G}_{2, N}(x, z'; \tau), \tilde{G}_{3, N}(x, z'; \tau))$,
we obtain (\ref{asymptotics of the gradient of the refracted wave}). 
Further, 
$$
G_{0}(x - \tilde{z}') = E_{0}(x - \tilde{z}')\frac{x - \tilde{z}'}{\v{x - \tilde{z}'}}.
$$ 
\end{Remark}

\par

\setcounter{equation}{0}
\section{Snell's law and the asymptotics of $\Phi_\tau$ and $\nabla_x\Phi_\tau$}
\label{Snell's law and the asymptotics of Phi_tau and nabla_xPhi_tau}

In this section, the Laplace integral 
(\ref{the Laplace integral for the refracted part}) is treated. First, we 
check the properties of function (\ref{the path length for Snell's law}), 
which describes Snell's law.

\begin{Lemma}\label{Snell's law}
The function $l_{x, y}$ in $\Bbb R^2$ defined by 
(\ref{the path length for Snell's law}) satisfies the following properties: 
\par\noindent
(1) Fix $x, y \in \Bbb R^3$ with $x_3 < 0$ and $y_3 > 0$. For the function
$l(x, y)$ defined by (\ref{definition of l(x, y)}), there exists
the unique point $z' \in \Bbb R^2$ satisfying $l(x, y) = l_{x, y}(z')$.
This point $z'$ is denoted by $z'(x, y)$. 
This point $z'(x, y)$ is on the line segment $x'y'$.
\par\noindent
(2) There exists a constant $C >0$ such that
$$
\sum_{i, j = 1}^{2}\frac{{\partial}^2l_{x, y}}{\partial{z_i}\partial{z_j}}(z'(x, y))
\xi_i\xi_j 
\geq C\v{\xi'}^2
\qquad(\xi' \in \Bbb R^2, x \in \overline{D}, y \in \overline{B}).
$$
\par\noindent
(3) The point $z'(x, y)$ is $C^\infty$ for 
$x, y \in \Bbb R^3$ with $x_3 < 0$ and $y_3 > 0$.
\end{Lemma}
Proof: In the beginning, we show 
\begin{align}
l_{x, y}
= \inf\{l_{x, y}(z') \vert z' \in \Bbb R^2, \v{z' - x'} \leq \v{x' - y'}, 
\v{z' - y'} \leq \v{x' - y'}\}.
\label{A set belonging maximal points}
\end{align} 
To obtain (\ref{A set belonging maximal points}), it suffices to show 
$l_{x, y}(z') > l_{x, y}(y') $ for  
$z' \in \Bbb R^2$ satisfying $\v{z' - x'} > \v{x' - y'}$, and
$l_{x, y}(z') > l_{x, y}(x') $ for 
$\v{z' - y'} > \v{x' - y'} $.
Suppose $\v{z' - x'} > \v{x' - y'} $, then it follows that
\begin{align*}
\v{x - \tilde{z}'} = \sqrt{\v{x' - z'}^2+x_3^2} > \sqrt{\v{x' - y'}^2+x_3^2}
= \v{x - \tilde{y}'}.
\end{align*}
Since $ \v{y - \tilde{z}'} = \sqrt{\v{y' - z'}^2+y_3^2} \geq \v{y_3} 
= \v{y - \tilde{y}'}$, we have
$$
l_{x, y}(z') = \frac{1}{\sqrt{\gamma_-}}\v{\tilde{z}'-x}
+\frac{1}{\sqrt{\gamma_+}}\v{\tilde{z}'-y}
> \frac{1}{\sqrt{\gamma_-}}\v{\tilde{y}'-x}
+\frac{1}{\sqrt{\gamma_+}}\v{\tilde{y}'-y}
= l_{x, y}(y').
$$
For $\v{z' - y'} > \v{x' - y'} $, 
$l_{x, y}(z') > l_{x, y}(x') $ is shown similarly, which yields
(\ref{A set belonging maximal points}). 

\par

From (\ref{A set belonging maximal points}) and $l_{x, y} \in C^{\infty}(\Bbb R^2)$, 
there exists a point $z' \in \Bbb R^2$ attaining the minimum $l(x, y)$ of $l_{x, y}$. 
Since this point $z'$ satisfies $\partial_{z'}l_{x, y}(z') = 0$, 
\begin{align}
\frac{\sin\theta_-}{\sqrt{\gamma_-}}\frac{z' - x'}{\v{{z}'-x'}}
+\frac{\sin\theta_+}{\sqrt{\gamma_+}}\frac{z' - y'}{\v{{z}'-y'}} = 0,
\label{Snell's law-1}
\end{align}
where $0 \leq \theta_{\pm} < \pi/2$ is taken by 
\begin{align}
\sin\theta_- = \frac{\v{z' - x'}}{\v{\tilde{z}'-x}}, \qquad
\sin\theta_+ = \frac{\v{z' - y'}}{\v{\tilde{z}'-y}}. \qquad
\label{Snell's law-2}
\end{align}
From (\ref{Snell's law-1}), it follows that $z'$ is on the segment $x'y'$, and
satisfies Snell's law, 
\begin{align*}
\frac{\sin\theta_-}{\sqrt{\gamma_-}} = \frac{\sin\theta_+}{\sqrt{\gamma_+}}.
\end{align*}
\par
We show that this point $z'$ is unique. If $x' = y'$, then $z' = x' = y'$, which
yields $\theta_\pm = 0$. Thus, $z'$ is uniquely determined. If $x' \neq y'$, 
this point is expressed by $z' = x'+t_0(y'-x')$ for some $0 < t_0 < 1$. 
We define $\varphi(t)$ by 
$\varphi(t) = l_{x, y}(x'+t(y'-x'))$.  Note that $t_0$ satisfies
$\varphi'(t_0) = 0$. Since $\varphi'(0) < 0$, $\varphi'(1) > 0$ and
$\varphi''(t) > 0$ for $0 \leq  t \leq 1$, there exists only one $0 < t_0 < 1$
with $\varphi'(t_0) = 0$, which yields uniqueness of $z'$.
Thus, we obtain (1) of Lemma \ref{Snell's law}.

\par

For (2), differentiate $l_{x, y}$ and obtain
\begin{align*}
\frac{{\partial}^2l_{x, y}}{\partial{z_i}\partial{z_j}}(z')
&= \frac{1}{\sqrt{\gamma_-}\v{\tilde{z}'-x}}
\Big\{\delta_{ij}-\Big(\frac{\v{z'-x'}}{\v{\tilde{z}'-x}}\Big)^2
\frac{(z_j - x_j)(z_i-x_i)}{\v{z'-x'}^2}\Big\}
\\&
\,\,\,
+
\frac{1}{\sqrt{\gamma_+}\v{\tilde{z}'-y}}
\Big\{\delta_{ij}-\Big(\frac{\v{z'-y'}}{\v{\tilde{z}'-y}}\Big)^2
\frac{(z_j - y_j)(z_i-y_i)}{\v{z'-y'}^2}\Big\}. 
\end{align*}
We put ${\bf e} = (z' - x')/\v{x' - z'}$. Then, (\ref{Snell's law-1}) and (\ref{Snell's law-2}) implies 
$(z' - y')/\v{y' - z'} = -{\bf e} $ since $z' = z(x, y)$ is in the 
line segment $x'y'$, which yields
\begin{align}
\displaystyle
\sum_{i, j = 1}^{2}\frac{{\partial}^2l_{x, y}}{\partial{z_i}\partial{z_j}}(z'(x, y))
\xi_i\xi_j 
&
=\frac{1}{\sqrt{\gamma_-}\v{\tilde{z}'-x}}
\Big(\v{\xi'}^2-(\sin\theta_-)^2({\bf e}\cdot\xi')^2\Big)
\label{form of Hessian of l}
\\
&
\,\,\,
+\frac{1}{\sqrt{\gamma_+}\v{\tilde{z}'-y}}
\Big(\v{\xi'}^2-(\sin\theta_+)^2({\bf e}\cdot\xi')^2\Big).
\nonumber
\end{align}
Hence, we have 
\begin{align}
\sum_{i, j = 1}^{2}&\frac{{\partial}^2l_{x, y}}{\partial{z_i}\partial{z_j}}(z'(x, y))
\xi_i\xi_j 
\geq 
\frac{(1-(\sin\theta_-)^2)\v{\xi'}^2
}{\sqrt{\gamma_-}\v{x_3}}
+
\frac{(1-(\sin\theta_+)^2)\v{\xi'}^2}{\sqrt{\gamma_+}\v{y_3}}.
\label{Snell's law-4}
\end{align}
\par

Since $\overline{D}\times\overline{B} \subset  \Bbb R^3_-\times\Bbb R^3_+$ is
bounded, there exist constants $L > 0$ and $A > 0$ such that
$$
\v{x' - y'} \leq L, A \leq \v{x_3} \leq A^{-1}, 
A \leq \v{y_3} \leq A^{-1}
\qquad
(x \in \overline{D}, y \in \overline{B}).
$$
Note that $t \mapsto \frac{t}{\sqrt{t^2+A^2}}$ is monotone increasing for
$t \geq 0$, for $z' = z'(x, y)$, it follows that
$$
0 \leq \sin\theta_\pm \leq \frac{L}{\sqrt{L^2+A^2}}, 
\quad(x \in \overline{D}, y \in \overline{B}).
$$
Combining the above estimate with (\ref{Snell's law-4}), we obtain
$$
\sum_{i, j = 1}^{2}
\frac{{\partial}^2l_{x, y}}{\partial{z_i}\partial{z_j}}
(z'(x, y))\xi_i\xi_j 
\geq \Big(\frac{1}{\sqrt{\gamma_-}}
+\frac{1}{\sqrt{\gamma_+}}\Big)\frac{A}{L^2+A^2}\v{\xi'}^2
\qquad(x \in \overline{D}, y \in \overline{B}, \xi' \in \Bbb R^2),
$$
which yields (2).

\par

Last, we show (3). We put $F(x, y, z') = l_{x, y}(z') $,
which is $C^\infty$ for $(x, y, z') \in \Bbb R^3_-\times\Bbb R^3_+\times\Bbb R^2$,
and $F(x, y, z'(x, y)) = 0$. Further, from (2), we obtain
${\rm det}\Big(\frac{\partial{F}}{\partial{z'}}
(z'(x, y))\Big) \neq 0$.
Hence implicit function theorem yields smoothness of $z'(x, y)$. 
\hfill$\blacksquare$
\vskip1pc\noindent

\par\noindent

Note that (\ref{form of Hessian of l}) implies 
\begin{align*}
\sum_{i, j = 1}^{2}\frac{{\partial}^2l_{x, y}}{\partial{z_i}\partial{z_j}}(z'(x, y))
\eta_i\xi_j 
&= \frac{1}{\sqrt{\gamma_-}\v{\tilde{z}'-x}}
\Big(\eta'\cdot\xi'-(\sin\theta_-)^2({\bf e}\cdot\eta')({\bf e}\cdot\xi')\Big)
\\
&
\,\,\,
+
\frac{1}{\sqrt{\gamma_+}\v{\tilde{z}'-y}}
\Big(\eta'\cdot\xi'-(\sin\theta_+)^2({\bf e}\cdot\eta')({\bf e}\cdot\xi')\Big)
\end{align*}
for $\xi' = (\xi_1, \xi_2)$ and $\eta' = (\eta_1, \eta_2)$. Hence, the eigenvalues 
and eigenvectors of
the Hessian $H(x, y) = {\rm Hess}(l_{x, y})(z'(x, y))$ are given by 
$$\left\{
\begin{array}{l}
\displaystyle
H(x, y)\frac{z' - x'}{\v{x' - z'}}= \Big(
\frac{x_3^2}{\sqrt{\gamma_-}\v{\tilde{z}'-x}^3}
+
\frac{y_3^2}{\sqrt{\gamma_+}\v{\tilde{z}'-y}^3}\Big)
\frac{z' - x'}{\v{x' - z'}},
\\
\\
\displaystyle
H(x, y){\bf e}'= \Big(
\frac{1}{\sqrt{\gamma_-}\v{\tilde{z}'-x}}
+ \frac{1}{\sqrt{\gamma_+}\v{\tilde{z}'-y}}\Big)
{\bf e}',
\end{array}
\right.
$$
where ${\bf e}'$ is a unit vector with ${\bf e}'\cdot(x'-y')/\v{x' - y'} = 0$.
Thus, we also obtain
\begin{align*}
{\rm det}H(x, y)
= \Big(\frac{1}{\sqrt{\gamma_-}\v{\tilde{z}'-x}}
+ \frac{1}{\sqrt{\gamma_+}\v{\tilde{z}'-y}}\Big)
\Big(\frac{x_3^2}{\sqrt{\gamma_-}\v{\tilde{z}'-x}^3}
+\frac{y_3^2}{\sqrt{\gamma_+}\v{\tilde{z}'-y}^3}\Big).
\end{align*}

\par

Now we are in the position to show 
Proposition \ref{Asymptotics of the refracted part of the gradient of the FS}.
Here, we need to assume $\gamma_+ < \gamma_-$. 
\par\noindent
Proof of Proposition \ref{Asymptotics of the refracted part of the gradient of the FS}.
From (\ref{expression of the FS for x_3 < 0}) and 
(\ref{asymptotics of the refracted wave}), 
we have
\begin{align}
\Phi_\tau(x, y)
&= \frac{{\tau}}{(4\pi)^2\gamma_-\gamma_+}\sum_{j = 0}^{N-1}\tau^{-j}
\int_{\Bbb R^2}e^{-{\tau}l_{x, y}(z')}f_{j}(z'; x, y)dz'
\label{the kernel representation of the FS}
\\
&
\,\,\,
+\frac{{\tau}}{(4\pi)^2\gamma_-\gamma_+}
\int_{\Bbb R^2}\frac{e^{-{\tau}l_{x, y}(z')}}{\v{x - \tilde{z}'}\v{\tilde{z}'-y}}
\tilde{E}_{N}(x, z'; \tau)dz', 
\nonumber
\end{align}
where 
\begin{align*}
f_{j}(z'; x, y) &= 
\frac{\gamma_-^{j/2}}{\v{x - \tilde{z}'}^{j+1}\v{\tilde{z}'-y}}
E_{j}(x - \tilde{z}')
\quad
(j = 0, 1, \ldots).
\end{align*}
From Proposition \ref{Asymptotics of refracted parts for case 1}, 
$E_{j}(x - \tilde{z}')$ are $C^\infty$ for $ x \in \Bbb R^3_-$ and
$z' \in \Bbb R^2$, which yields $f_{j} \in C^\infty(\Bbb R^2\times\Bbb R^3_-\times\Bbb R^3_+)$. 
Hence, we can apply (\ref{asymptotic formula of I(tau; x, y)})
for each Laplace integral with $f_j$, and obtain the asymptotics. 
For the remainder term of (\ref{the kernel representation of the FS}), 
\begin{align*}
\left\vert
\int_{\Bbb R^2}\frac{e^{-{\tau}l_{x, y}(z')}}{\v{x - \tilde{z}'}\v{\tilde{z}'-y}}
\tilde{E}_{N}(x, z',; \tau)dz'
\right\vert
&\leq C_N\frac{\gamma_-^{N/2}e^{-{\tau}l(x, y)}}{\tau^N}
\int_{\Bbb R^2}\frac{dz'}{\v{x - \tilde{z}'}^{N+1}\v{\tilde{z}'-y}}
\\&
\leq C_N\frac{\gamma_-^{N/2}e^{-{\tau}l(x, y)}}{\tau^N}
\qquad(x \in \overline{D}, y \in \overline{B}, \tau \geq 1). 
\end{align*}
Thus, we obtain the asymptotics for $\Phi_\tau(x, y) $ in 
Proposition \ref{Asymptotics of the refracted part of the gradient of the FS} 
uniformly in $x \in \overline{D}, y \in \overline{B}$,
where $\Phi_{j}^{(0)}(x, y)$ is given by
$\Phi_{j}^{(0)}(x, y) = \sum_{p = 0}^{j}(L_pf_{j - p}(\cdot; x, y))(z'(x, y))$.
Hence, we also have 
$\Phi_{0}^{(0)}(x, y) = f_0(z'(x, y); x, y)$, 
which yields (\ref{the form of Phi_{0}^{(0)}}).

\par
For $\nabla_x\Phi_\tau(x, y)$, 
differentiating (\ref{expression of the FS for x_3 < 0}), and using 
(\ref{asymptotics of the gradient of the refracted wave}), we obtain
\begin{align*}
\nabla_x\Phi_\tau(x, y)
&= \frac{-{\tau}^2}{(4\pi)^2\gamma_-^{3/2}\gamma_+}\sum_{j = 0}^{N-1}\tau^{-j}
\int_{\Bbb R^2}e^{-{\tau}l_{x, y}(z')}g_{j}(z'; x, y)dz'
\\&
\,\,\,
-\frac{{\tau}^2}{(4\pi)^2\gamma_-^{3/2}\gamma_+}
\int_{\Bbb R^2}\frac{e^{-{\tau}l_{x, y}(z')}}{\v{x - \tilde{z}'}\v{\tilde{z}'-y}}
\tilde{G}_{N}(x, z'; \tau)dz', 
\nonumber
\end{align*}
where 
\begin{align*}
g_{j}(z'; x, y) &= 
\frac{\gamma_-^{j/2}}{\v{x - \tilde{z}'}^{j+1}\v{\tilde{z}'-y}}
G_{j}(x - \tilde{z}')
\quad
(j = 0, 1, \ldots).
\end{align*}
From Proposition \ref{Asymptotics of refracted parts for case 1}, 
$G_{j}(x - \tilde{z}')$ are $C^\infty$ for $ x \in \Bbb R^3_-$ and
$z' \in \Bbb R^2$, which yields $g_{j} \in C^\infty(\Bbb R^2\times\Bbb R^3_-\times\Bbb R^3_+)$. 
Hence, as for $\Phi_\tau(x, y)$, we obtain the asymptotics for
$\nabla_x\Phi_\tau(x, y)$ described in 
Proposition \ref{Asymptotics of the refracted part of the gradient of the FS}. 
In this case, $\Phi_{j}^{(1)}(x, y)$ is given by
$\Phi_{j}^{(1)}(x, y) = \sum_{p = 0}^{j}(L_pg_{j - p}(\cdot; x, y))(z'(x, y))$.
From this and Remark \ref{Asymptotics of gradient of refracted parts for case 1},
we obtain 
$$
\Phi_{0}^{(1)}(x, y) = g_{0}(z'(x, y); x, y) 
=  E_{0}(x - \tilde{z}')\frac{x - \tilde{z}'}{\v{x - \tilde{z}'}^2\v{y - \tilde{z}'}}, 
$$
which yields (\ref{the form of Phi_{1}^{(0)}}). 
This completes the proof of 
Proposition \ref{Asymptotics of the refracted part of the gradient of the FS}.
\hfill$\blacksquare$
\vskip1pc\noindent

\setcounter{equation}{0}
\section{Proof of Theorem \ref{Theorem 1.2}}
\label{Proof of Theorem ref{Theorem 1.2}}

We put $L(x, y, \xi) = l(x, y)+l(x, \xi)$, 
\begin{align*}
\tilde{K}(x, y, \xi, \tau) &= {\tau}\big(\Phi_{0}^{(1)}(x, y){\cdot}Q_{0, \tau}^{(1)}(x, \xi)
+Q_{0, \tau}^{(1)}(x, y){\cdot}\Phi_{0}^{(1)}(x, \xi)
+Q_{0, \tau}^{(1)}(x, y){\cdot}Q_{0, \tau}^{(1)}(x, \xi)\big)
\end{align*} 
$$
J(x, y, \xi) = \frac{(8\pi)^{-2}\gamma_+^{-2}(\gamma_-)^{-3}}{
\sqrt{{\rm det}H(x,y){\rm det}H(x, \xi)}},
\quad
h_0(x, y, \xi) = 
\frac{x - \tilde{z}'(x, y)}{\v{x - \tilde{z}'(x, y)}} 
\cdot\frac{x - \tilde{z}'(x, \xi)}{\v{x - \tilde{z}'(x, \xi)}}
$$
and $h_1(x, y, \xi) = \Phi_{0}^{(0)}(x, y)\Phi_{0}^{(0)}(x, \xi)$.
From Proposition \ref{Asymptotics of the refracted part of the gradient of the FS}
for $N = 0$, it follows that
\begin{align}
\nabla_x\Phi_\tau(x, y)\cdot\nabla_x\Phi_\tau(x, \xi) 
= {\tau}^2{e}^{-{\tau}L(x, y, \xi)}&J(x, y, \xi)K(x, y, \xi, \tau)
\label{the form of the gradient^2 of FS}
\\
&\quad\quad\quad
(x \in \Bbb R^3_-, y, \xi \in \Bbb R^3_+), 
\nonumber
\end{align}
where $K(x, y, \xi, \tau) = h_0(x, y, \xi)h_1(x, y, \xi)+\tau^{-1}\tilde{K}(x, y, \xi, \tau)$.
From (\ref{the form of E_0}), there exist constants $0 < C_1 < C_2$ such that
\begin{align}
C_1 \leq J(x, y, \xi) &\leq C_2 
\qquad(x \in \overline{D}, y, \xi \in \overline{B}), 
\label{Proof of asymptotics for the norm of v 1}
\\
C_1 \leq h_1(x, y, \xi) &\leq C_2 
\qquad(x \in \overline{D}, y, \xi \in \overline{B}), 
\label{Proof of asymptotics for the norm of v 3}
\\
\v{\tilde{K}(x, y, \xi, \tau)} &\leq C_2 
\qquad(x \in \overline{D}, y, \xi \in \overline{B}, \tau > 1).
\label{Proof of asymptotics for the norm of v 2}
\end{align}

\par\noindent
\begin{Lemma}\label{properties of the points giving the minimum}
Put $l_0 = \min_{x \in \overline{D}, y \in \overline{B}}l(x, y) > 0$ and
$l_1 = \min_{y, \xi \in \overline{B}, x \in \overline{D}}L(x, y, \xi)$.
\par\noindent
(1) $ l_1 = 2l_0$. Further, if $x_0 \in \overline{D}$, $y_0 \in \overline{B}$
satisfy $l_0 = l(x_0, y_0)$, then $L(x_0, y_0, y_0) = l_1$.
\par\noindent
(2) If $x_0 \in \overline{D}$, $y_0 \in \overline{B}$ satisfy $l_0 = l(x_0, y_0)$,
then $x_0 \in \partial{D}$ and $y_0 \in \partial{B}$.
Further, 
$$
\nu_{x_0} = \frac{\tilde{z}'(x_0, y_0)-x_0}{\v{x_0 - \tilde{z}'(x_0, y_0)}},
\qquad
\nu_{y_0} = \frac{\tilde{z}'(x_0, y_0)-y_0}{\v{y_0 - \tilde{z}'(x_0, y_0)}},
$$
where $\nu_{x_0}$ and $\nu_{y_0}$ are the unit outer normal of
$\partial{D}$ and $\partial{B}$ at $x_0$ and $y_0$, respectively.
\par\noindent
(3) If $x_1 \in \overline{D}$, $y_1, \xi_1 \in \overline{B}$ satisfy
$L(x_1, y_1, \xi_1) = l_1$, then $y_1 = \xi_1$ and 
$l(x_1, y_1) = l_0 $.
\end{Lemma}
Proof: For any $y, \tilde{y} \in \overline{B}$ and $x \in \overline{D} $,
$L(x, y, \tilde{y}) = l(x, y)+l(x, \tilde{y}) \geq 2l_0$. Since 
$l(x, y)$ is continuous on the compact set $\overline{D}\times\overline{B}$,
we can take points $x_0 \in \overline{D}$ and $y_0 \in \overline{B}$ satisfying
$l_0 = l(x_0, y_0)$, which implies $2l_0 \leq L(x_0, y_0, y_0) = 2l_0$.
Thus, we obtain (1).
\par
To show (2), assume that $x_0 \in \overline{D}$, $y_0 \in \overline{B}$ satisfy $l_0 = l(x_0, y_0)$. 
If $x_0 \notin \partial{D}$, there exists $\delta > 0$ such that
$B_{2\delta}(x_0) \subset D$, where for $a \in \Bbb R^3$ and $r > 0$, we put
$B_r(a) = \{ x \in \Bbb R^3 \vert \v{x - a} < r\}$.
We put $z'_0 = z'(x_0, y_0)$, $\tilde{z}_0' = (z_0, 0) \in \Bbb R^3$,
${\bf e} = (\tilde{z}_0'-x_0)/\v{\tilde{z}_0'-x_0}$ and
$x_1 = x_0 + \delta{\bf e}$. Then, $x_1 \in D $ for $\v{x_1 - x_0} < \delta$, and
$\tilde{z}_0' - x_1 = \tilde{z}_0' - (x_0 + \delta{\bf e}) = \v{\tilde{z}_0'-x_0}{\bf e} - \delta{\bf e}
= (\v{\tilde{z}_0'-x_0} - \delta){\bf e}$. These imply
\begin{align*}
l_{x_1, y_0}(z'_0) &= \frac{1}{\sqrt{\gamma_-}}\v{\tilde{z}_0'-x_1}
+\frac{1}{\sqrt{\gamma_+}}\v{\tilde{z}_0'-y_0}
= \frac{1}{\sqrt{\gamma_-}}(\v{\tilde{z}_0'-x_0}-\delta)
+\frac{1}{\sqrt{\gamma_+}}\v{\tilde{z}_0'-y_0}
\\&
< l_{x_0, y_0}(z'_0) = l_0,
\end{align*}
which is contradiction. Thus, we obtain $x_0 \in \partial{D}$. Similarly, we have
$y_0 \in \partial{B}$.
\par
Next, we show ${\bf e}$ is a unit outer normal of $\partial{D}$ at $x_0$.
Take any $C^1$ class curve 
$c: (-\varepsilon, \varepsilon) \to \partial{D}$ with $c(0) = x_0$.
Since 
$$
l_0 \leq l(c(t), y_0) 
\leq l_{c(t), y_0}(z_0') = 
\frac{1}{\sqrt{\gamma_-}}\v{\tilde{z}_0'-c(t)}
+ \frac{1}{\sqrt{\gamma_+}}\v{\tilde{z}_0'-y_0}
$$
and $l_{c(0), y_0}(z_0') = l_{x_0, y_0}(z_0') = l_0$,
the function $(-\varepsilon, \varepsilon) \ni t \mapsto l_{c(t), y_0}(z_0')$
take a minimum at $t = 0$. This implies
$$
0 = \frac{d}{dt}(l_{c(t), y_0}(z_0'))\big\vert_{t = 0}
= \frac{1}{\sqrt{\gamma_-}}
\frac{c(0) - \tilde{z}_0'}{\v{\tilde{z}_0'-c(0)}}{\cdot}c'(0)
= \frac{1}{\sqrt{\gamma_-}}{\bf e}{\cdot}c'(0),
$$
which yields ${\bf e}$ is a unit normal of $\partial{D}$.

\par
To obtain ${\bf e}$ is outward, it suffices to show $x_0+\delta{\bf e} \notin D$
for $\delta > 0$ small enough. For any $0 < \delta < \v{x_0-\tilde{z}_0'}$, it follows that
\begin{align*}
l_{x_0+\delta{\bf e}, y_0}(\tilde{z}_0')
&= \frac{1}{\sqrt{\gamma_-}}\v{x_0+\delta{\bf e}-\tilde{z}_0'}
+ \frac{1}{\sqrt{\gamma_+}}\v{\tilde{z}_0'-y_0}
\\&
= \frac{1}{\sqrt{\gamma_-}}(\v{x_0-\tilde{z}_0'}-\delta)
+ \frac{1}{\sqrt{\gamma_+}}\v{\tilde{z}_0'-y_0}
< l_0,
\end{align*}
which yields $x_0+\delta{\bf e} \notin D$. For $y_0 \in \partial{B}$, we can show
similarly, which obtain (2). 
\par

Last, we show (3). 
Take $x_1 \in \overline{D}$ and $y_1, \xi_1 \in \overline{B}$ with
$L(x_1, y_1, \xi_1) = l_1$. From 
$l_0 \leq l(x_1, y_1)$, $l_0 \leq l(x_1, \xi_1)$, it follows that
$$
l_0 \leq l(x_1, y_1) = L(x_1, y_1, \xi_1) - l(x_1, \xi_1)
= 2l_0 - l(x_1, \xi_1) \leq 2l_0 - l_0 = l_0,
$$
which yields $l(x_1, y_1) = l_0$. We can obtain
$l(x_1, \xi_1) = l_0$ similarly. To finish the proof, it 
suffices to show $y_1 = \xi_1$.
\par
We put $z_1' = z'(x_1, y_1) \in \Bbb R^2$, $\tilde{z}_1' = (z_1', 0)$, 
$\eta_1' = z'(x_1, \xi_1) \in \Bbb R^2$ and $\tilde{\eta}_1' = (\eta_1', 0)$, 
where $z'(x_1, y_1)$ and $z'(x_1, \xi_1) \in \Bbb R^2$ are determined by (1) of 
Lemma \ref{Snell's law}, respectively. 
From (2) of Lemma \ref{properties of the points giving the minimum}, 
$\nu_{x_1} = (\tilde{z}_1'-x_1)/\v{\tilde{z}_1'-x_1} = (\tilde{\eta}_1'-x_1)/\v{\tilde{\eta}_1'-x_1}$.
Taking the inner product of this vector and $(0, 0, 1) \in \Bbb R^3$, we obtain
$ x_{1, 3}/\v{\tilde{z}_1'-x_1} = x_{1,3}/\v{\tilde{\eta}_1'-x_1}$, where
$x_1 = (x_{1, 1}, x_{1, 2}, x_{1, 3})$. 
Since $x_{1, 3} \neq 0$, it follows that $\v{\tilde{z_1'}-x_1} = \v{\tilde{\eta_1'}-x_1}$,
which yields $\tilde{z}_1' = \tilde{\eta}_1'$. From this and $l(x_1, y_1) = l(x_1, \xi_1) = l_0$,
$\v{\tilde{z}_1'-y_1} = \v{\tilde{\eta}_1'-\xi_1}$ also follows. 
\par
Now, we remember Snell's law,  
$$
\frac{1}{\sqrt{\gamma_-}}\frac{z_1' - x_1'}{\v{\tilde{z}_1'-x_1}}
+\frac{1}{\sqrt{\gamma_+}}\frac{z_1' - y_1'}{\v{\tilde{z}_1'-y_1}} = 0,
\qquad
\frac{1}{\sqrt{\gamma_-}}\frac{\eta_1' - x_1'}{\v{\tilde{\eta}_1'-x_1'}}
+\frac{1}{\sqrt{\gamma_+}}\frac{\eta_1' - \xi_1'}{\v{\tilde{\eta}_1'-\xi_1}} = 0,
$$
which are derived from (\ref{Snell's law-1}) and (\ref{Snell's law-2}).
These relations imply $y_1' = \xi_1'$. Since $\v{\tilde{z}_1'-y_1} = \v{\tilde{\eta}_1'-\xi_1}$,
it follows that $y_{1, 3}^2 = \xi_{1, 3}^2$, which yields $y_{1, 3} = \xi_{1, 3}$
since both are positive. Thus, we obtain $y_1 = \xi_1$, which 
completes the proof of Lemma \ref{properties of the points giving the minimum}.
\hfill$\blacksquare$
\vskip1pc

Now, we are in a position to show Theorem \ref{Theorem 1.2}. Combining 
(\ref{the form of the gradient^2 of FS})-(\ref{Proof of asymptotics for the norm of v 2})
and (\ref{the integral of the gradient of v in D}) 
with (1) of Lemma \ref{properties of the points giving the minimum}, we obtain the estimate
of the right side in (\ref{1.7}). The problem is
to show the left side of (\ref{1.7}).
\par
Since 
$h_0(x, y, y) = 1$ for $x \in \Bbb R^3_-$, $y \in \Bbb R^3_+$, and
$h_0(x, y, \xi)$ is continuous for $x \in \Bbb R^3_-$, $y, \xi \in \Bbb R^3_+$,
which is from Lemma \ref{Snell's law}, there exists a constant $\delta > 0$ 
such that 
\begin{equation}
\text{$h_0(x, y, \xi) \geq 1/2$ for $\v{y - \xi} < 3\delta$, $y, \xi \in \overline{B}$, $x \in \overline{D}$. }
\label{(i) for the main estimate}
\end{equation}
We can take $\delta > 0$ in (\ref{(i) for the main estimate}) sufficiently small 
to be $\overline{B_{4\delta}(y)} \subset \Bbb R^3_+$ $(y \in \overline{B})$.
\par

Put
$E = \{(x_0, y_0, \xi_0) \in \overline{D}\times\overline{B}\times\overline{B} \vert 
L(x_0, y_0, \xi_0) = 2l_0\}$ and 
$E_0 = \{(x_0, y_0) \in \overline{D}\times\overline{B} \vert l(x_0, y_0) = l_0\}$.
(1) and (3) of Lemma \ref{properties of the points giving the minimum} imply 
$E \neq \emptyset$ and $E_0 \neq \emptyset$, and
$E = \{(x_0, y_0, \xi_0) \in \overline{D}\times\overline{B}\times\overline{B} \vert 
l(x_0, y_0) = l_0, y_0 = \xi_0 \}$. 
Since $E_0$ is a compact set, there exist finite points
$(x_0^{(k)}, y_0^{(k)}) \in E_0$ $(k = 1, 2, \ldots, N)$ such that 
$E_0 \subset \cup_{k = 1}^{N}
B_\delta(x_0^{(k)}){\times}B_\delta(y_0^{(k)})$, where $\delta > 0$ is given in
(\ref{(i) for the main estimate}). Then, $E \subset \cup_{k = 1}^{N}
B_\delta(x_0^{(k)}){\times}B_\delta(y_0^{(k)}){\times}B_\delta(y_0^{(k)})$ holds.
\par
We put 
${\mathcal W} = \cup_{k = 1}^{N}
B_\delta(x_0^{(k)}){\times}B_\delta(y_0^{(k)}){\times}B_\delta(y_0^{(k)})
\subset \Bbb R^3_+{\times}\Bbb R^3_-{\times}\Bbb R^3_-$.
Since $\v{y - \xi} \leq \v{y - y_0^{(k)}}+\v{y_0^{(k)}-\xi} < 2\delta$
for $(y, \xi) \in B_\delta(y_0^{(k)}){\times}B_\delta(y_0^{(k)})$, 
from (\ref{(i) for the main estimate}), it follows that  
\begin{equation}
\text{$h_0(x, y, \xi) \geq 1/2$ 
$((x, y, \xi) \in {\mathcal W})$}
\label{(ii) for the main estimate}
\end{equation}
Since $E \subset {\mathcal W}$, $L(x, y, \xi) = l(x, y)+l(x, \xi) > l_1 = 2l_0$
on the compact set $\overline{D}{\times}\overline{B}{\times}\overline{B}
\setminus{\mathcal W}$. Thus, there exists  $c_0 > 0$ such that
\begin{equation}
L(x, y, \xi) \geq 2l_0+c_0  \qquad
(x, y, \xi) \in \overline{D}{\times}\overline{B}{\times}\overline{B}
\setminus{\mathcal W}.
\label{(iii) for the main estimate}
\end{equation}
From (\ref{(ii) for the main estimate}), (\ref{Proof of asymptotics for the norm of v 3}) and
(\ref{Proof of asymptotics for the norm of v 2}), 
for any $(x, y, \xi) \in {\mathcal W}$, it follow that
\begin{align*}
K(x, y, \xi, \tau) &= h_0(x, y, \xi)h_1(x, y, \xi)
+\tau^{-1}\tilde{K}(x, y, \xi, \tau)
\geq \frac{C_1}{2}-\tau^{-1}C_2,
\end{align*}
which yields that there exists a constant $\tau_0 \geq 1$ such that 
$K(x, y, \xi, \tau) \geq C_1/4$ $((x, y, \xi) \in {\mathcal W}, 
\tau \geq \tau_0)$.
From (\ref{Proof of asymptotics for the norm of v 3}) and
(\ref{Proof of asymptotics for the norm of v 2}), for any
$(x, y, \xi) \in \overline{D}\times\overline{B}\times\overline{B}$,
it follows that
\begin{align*}
\v{K(x, y, \xi, \tau)} &\leq \v{h_0(x, y, \xi)}\v{h_1(x, y, \xi)}
+\tau^{-1}\v{\tilde{K}(x, y, \xi, \tau)}
\leq C_2(1+\tau^{-1}),
\end{align*}
which imples $\v{K(x, y, \xi, \tau)} \leq 2C_2$ 
$((x, y, \xi) \in \overline{D}\times\overline{B}\times\overline{B}, 
\tau \geq 1)$.

\par

From the above estimates of $K$, and 
(\ref{the form of the gradient^2 of FS}), 
(\ref{(iii) for the main estimate}), 
(\ref{Proof of asymptotics for the norm of v 1}), 
(\ref{Proof of asymptotics for the norm of v 2}), 
(\ref{the integral of the gradient of v in D}) and 
the assumption for $f$, it follows that
\begin{align*}
\int_{D}\v{\nabla_xv(x)}^2dx &\geq {\tau}^2C_0^2C_1\frac{C_1}{4}\int_{{\mathcal W}}
{e}^{-{\tau}L(x, y, \xi)}dyd{\xi}dx
\\&
\,\,\,
- {\tau}^2e^{-{\tau}(2l_0+c_0)}2C_2^2
\int_{B{\times}B}dyd{\xi}\v{f(y)}\v{f(\xi)}
\int_{D}dx
\\&
\geq
{\tau}^2\Big(\frac{C_0^2C_1^2}{4}\int_{{\mathcal W}}
{e}^{-{\tau}L(x, y, \xi)}dyd{\xi}dx 
- C_3e^{-{\tau}(2l_0+c_0)}\Big)
\quad(\tau \geq \tau_0),
\end{align*}
where $C_3 = 2C_2^2{\rm Vol(D)}{\rm Vol(B)}\V{f}_{L^2(B)}^2$.
Thus, to obtain the left hand of (\ref{1.7}), it suffices to show the following estimate:
\begin{Lemma}\label{estimates of the integral as the main term}
There exist $\tau_1 \geq 1$ and $C > 0$ such that 
$$
\int_{{\mathcal W}}
{e}^{-{\tau}L(x, y, \xi)}dyd{\xi}dx \geq C\tau^{-6}e^{-2l_0\tau}
\qquad(\tau \geq \tau_1).
$$
\end{Lemma}
Proof: In what follows, we write $x_0^{(1)}$ and $y_0^{(1)}$ as $x_0$ and $y_0$,
respectively. Since 
(2) of Lemma \ref{properties of the points giving the minimum} implies
$x_0 \in \partial{D}$ and $y_0 \in \partial{B}$, and $\partial{D}$ and
$\partial{B}$ are $C^1$ surfaces, there exist $p_0$, $q_0 \in \Bbb R^3$ and 
$r > 0$ such that $B_r(p_0) \subset B_\delta(x_0)$,
$B_r(q_0) \subset B_\delta(y_0)$, $x_0 \in \partial{B_r(p_0)}$ and 
$y_0 \in \partial{B_r(q_0)}$. Then, 
${\mathcal W} \supset 
B_\delta(x_0){\times}B_\delta(y_0){\times}B_{\delta(y_0)}
\supset {B_r(p_0)}{\times}{B_r(q_0)}{\times}{B_r(q_0)}$, 
and 
\begin{align*}
\int_{{\mathcal W}}{e}^{-{\tau}L(x, y, \xi)}dyd{\xi}dx
&\geq
\int_{{B_r(p_0)}{\times}{B_r(q_0)}{\times}{B_r(q_0)}}
{e}^{-{\tau}(l(x, y)+l(x, \xi))}dyd{\xi}dx
\\&
= 
\int_{B_r(p_0)}dx\Big(\int_{B_r(q_0)}{e}^{-{\tau}l(x, y)}dy\Big)^2.
\end{align*}
\par
We put $z_0' = {z}'(x_0, y_0) \in \Bbb R^2$ and 
$\tilde{z}_0 = (z_0', 0) \in \partial\Bbb R^3_+$.
Since
$$
l(x, y) \leq l_{x, y}(z_0')
= \frac{1}{\sqrt{\gamma_-}}\v{\tilde{z}_0-x}
+\frac{1}{\sqrt{\gamma_+}}\v{\tilde{z}_0-y},
$$
it follows that
\begin{align*}
\int_{{\mathcal W}}{e}^{-{\tau}L(x, y, \xi)}dyd{\xi}dx
&\geq 
\int_{B_r(p_0)}{e}^{-\frac{2{\tau}}{\sqrt{\gamma_-}}\v{\tilde{z}_0-x}}dx
\Big(\int_{B_r(q_0)}{e}^{-\frac{{\tau}}{\sqrt{\gamma_+}}\v{\tilde{z}_0-y}}dy\Big)^2.
\end{align*}
As in (2) of Lemma \ref{properties of the points giving the minimum},
$\nu_{y_0} = (\tilde{z}_0-y_0)/\v{\tilde{z}_0-y_0}$, which yields 
$\inf_{y \in B_r(q_0)}\v{\tilde{z}_0 - y} = \v{\tilde{z}_0-y_0}$. Then, from 
Proposition 3.2 of \cite{Ikehata-Kawashita1}, it follows that
there exist constants $C > 0$ and $\tau_2 > 0$ such that 
$$
\int_{B_r(q_0)}{e}^{-\frac{{\tau}}{\sqrt{\gamma_+}}\v{\tilde{z}_0-y}}dy
\geq C\tau^{-2}{e}^{-\frac{{\tau}}{\sqrt{\gamma_+}}\v{\tilde{z}_0-y_0}}
\qquad(\tau \geq \tau_2).
$$
Similarly, taking the constants $C$ and $\tau_2$ larger if necessary, 
we also obtain
$$
\int_{B_r(p_0)}{e}^{-\frac{2{\tau}}{\sqrt{\gamma_-}}\v{\tilde{z}_0-x}}dx
\geq C\tau^{-2}{e}^{-\frac{2{\tau}}{\sqrt{\gamma_-}}\v{\tilde{z}_0-x_0}}
\qquad(\tau \geq \tau_2).
$$
Hence, we have
\begin{align*}
\int_{{\mathcal W}}&{e}^{-{\tau}L(x, y, \xi)}dyd{\xi}dx
\geq 
C^3\tau^{-6}{e}^{-\frac{2{\tau}}{\sqrt{\gamma_-}}\v{\tilde{z}_0-x_0}}
\Big({e}^{-\frac{{\tau}}{\sqrt{\gamma_+}}\v{\tilde{z}_0-y_0}}\Big)^2
= C^3\tau^{-6}{e}^{-2{\tau}l_0} 
\end{align*}
since
$$
\frac{\v{\tilde{z}_0-x_0}}{\sqrt{\gamma_-}}+
\frac{\v{\tilde{z}_0-y_0}}{\sqrt{\gamma_+}}
= l_{x_0, y_0}(z_0') = l(x_0, y_0) = l_0,
$$
which completes the proof of
Lemma \ref{estimates of the integral as the main term}.
\hfill
$\blacksquare$

\setcounter{equation}{0}

\vskip2pc

\centerline{{\bf Acknowledgements}}

MI was partially supported by JSPS KAKENHI Grant Number JP17K05331.
MK was partially supported by JSPS KAKENHI Grant Number JP16K05232.

\setcounter{equation}{0}
\appendix
\renewcommand{\theequation}{A.\arabic{equation}}

\section{Appendix. Proof of Lemma \ref{Lemma 1.1}}

The function $w$ defined by (\ref{definition of w}) satisfies
\begin{equation}
\begin{array}{ll}
\displaystyle
(\nabla\cdot\gamma\nabla-\tau^2)w+f=e^{-\tau T}F
& 
\displaystyle
\text{in}\,\Bbb R^3,
\end{array}
\label{A.1}
\end{equation}
where
\begin{equation}
\displaystyle
F=F(x,\tau)
=\partial_tu(x,T)+\tau u(x,T).
\label{A.2}
\end{equation}
Define
$$\displaystyle
R=w-v.
$$
It is easy to derive the following decomposition formula
of the indicator function which formally corresponds to the case when 
$\Omega=\Bbb R^3$  on 
(3.2) of Proposition 3.1 in \cite{IE}.
\begin{Prop}\label{Proposition A.1}
We have
\begin{equation}
\begin{array}{ll}
\displaystyle
\int_{\Bbb R^3}f Rdx &
\displaystyle
=\int_{\Bbb R^3}(\gamma_0I_3-\gamma)\nabla v\cdot\nabla v dx
+\int_{\Bbb R^3}\gamma\nabla R\cdot\nabla R dx
+\tau^2\int_{\Bbb R^3}\vert R\vert^2 dx
\\
\displaystyle
&
\displaystyle
\,\,\,
+
e^{-\tau T}\left(\int_{\Bbb R^3}FRdx-\int_{\Bbb R^3}Fvdx\right).
\end{array}
\label{A.3}
\end{equation}
\end{Prop}

It follows from (\ref{A.1}) that $w$ satisfies
$$\begin{array}{ll}
(\nabla\cdot\gamma\nabla-\tau^2)w
+\tilde{f}=0 &
\displaystyle
\text{in}\,\Bbb R^3,
\end{array}
$$
where
$$\displaystyle
\tilde{f}=f-e^{-\tau T}F.
$$
And also $v$ satisfies
$$\begin{array}{ll}
\displaystyle
(\nabla\cdot\gamma_0\nabla-\tau^2)v+\tilde{f}=e^{-\tau T}\tilde{F}
&
\displaystyle
\text{in}\,\Bbb R^3,
\end{array}
$$
where
$$\displaystyle
\tilde{F}=-F.
$$
Thus, changing the role of $v$ and $w$ in (\ref{A.3}), we obtain
\begin{align*}
\int_{\Bbb R^3}\tilde{f}(-R)dx &
\displaystyle
=-\int_{\Bbb R^3}(\gamma_0I_3-\gamma)\nabla w\cdot\nabla w dx
+\int_{\Bbb R^3}\gamma_0\nabla R\cdot\nabla R dx
+\tau^2\int_{\Bbb R^3}\vert R\vert^2 dx
\\&
\,\,\,\hskip5mm
+
e^{-\tau T}\left(\int_{\Bbb R^3}\tilde{F}(-R)dx-\int_{\Bbb R^3}\tilde{F}wdx\right).
\end{align*}
This is noting but the following formula.

\begin{Prop}\label{Proposition A.2}
We have
\begin{equation}
\begin{array}{ll}
\displaystyle
-\int_{\Bbb R^3}f Rdx &
\displaystyle
=\int_{\Bbb R^3}(\gamma-\gamma_0I_3)\nabla w\cdot\nabla w dx
+\int_{\Bbb R^3}\gamma_0\nabla R\cdot\nabla R dx
+\tau^2\int_{\Bbb R^3}\vert R\vert^2 dx
\\[5mm]
\displaystyle
&
\displaystyle
\,\,\,
+
e^{-\tau T}\left(\int_{\Bbb R^3}FRdx+\int_{\Bbb R^3}Fvdx\right).
\end{array}
\label{A.4}
\end{equation}
\end{Prop}

Now we are ready to prove (\ref{1.5}) and (\ref{1.6}).
It follows from (\ref{1.2}) that
$$\displaystyle
\int_{\Bbb R^3}(\gamma_0\nabla v\cdot\nabla v+\tau^2v^2-fv)dx=0,
$$
that is
$$\displaystyle
\int_{\Bbb R^3}\left\{\gamma_0\nabla v\cdot\nabla v+\left(\tau v-\frac{f}{2\tau}\right)^2\right\}dx
=\frac{1}{4\tau^2}\int_{\Bbb R^3}f^2dx.
$$
This yields, as $\tau\longrightarrow\infty$
$$\displaystyle
\int_{\Bbb R^3}(\gamma_0\nabla v\cdot\nabla v+\tau^2 v^2)dx=O(\tau^{-2}).
$$
Hence we obtain, as $\tau\longrightarrow\infty$
\begin{equation}
\displaystyle
\Vert v\Vert_{L^2(\Bbb R^3)}=O(\tau^{-2}).
\label{A.5}
\end{equation}
and
\begin{equation}
\displaystyle
\Vert \nabla v\Vert_{L^2(\Bbb R^3)}=O(\tau^{-1}).
\label{A.6}
\end{equation}
Rewrite (\ref{A.3}) as
\begin{equation}
\begin{array}{l}
\displaystyle
\int_{\Bbb R^3}\gamma\nabla R\cdot\nabla R dx
+\int_{\Bbb R^3}\left(\tau R-\frac{f-e^{-\tau t}F}{2\tau}\right)^2dx
\\[5mm]
\displaystyle
=-\int_{\Bbb R^3}(\gamma_0I_3-\gamma)\nabla v\cdot\nabla v dx
+e^{-\tau T}\int_{\Bbb R^3}Fvdx+\frac{1}{4\tau^2}\int_{\Bbb R^3}(f-e^{-\tau T}F)^2dx.
\end{array}
\label{A.7}
\end{equation}
Here from (\ref{A.2}) we have 
\begin{equation}
\Vert F\Vert_{L^2(\Bbb R^3)}=O(\tau).
\label{A.8}
\end{equation}
This together with (\ref{A.5}) and (\ref{A.6})
yields that the right-hand side on (\ref{A.7}) has a bound $O(\tau^{-2})$.  Thus we obtain
$$
\displaystyle
\int_{\Bbb R^3}\left(\gamma\nabla R\cdot\nabla R+\tau^2 R^2\right)\,dx=
O(\tau^{-2})
$$
and, in particular,
\begin{equation}
\displaystyle
\Vert R\Vert_{L^2(\Bbb R^3)}=O(\tau^{-2}).
\label{A.9}
\end{equation}
Now applying (\ref{A.5}) , (\ref{A.8})  and (\ref{A.9}) to the fourth-term in the right-hand side on 
(\ref{A.3}) we obtain 
\begin{equation}
\begin{array}{ll}
\displaystyle
\int_{\Bbb R^3}f Rdx &
\displaystyle
=\int_{\Bbb R^3}(\gamma_0I_3-\gamma)\nabla v\cdot\nabla v dx
\\
\\
\displaystyle
&
\,\,\,
\displaystyle
+\int_{\Bbb R^3}\gamma\nabla R\cdot\nabla R dx
+\tau^2\int_{\Bbb R^3}\vert R\vert^2 dx
+O(\tau^{-1}e^{-\tau T}).
\end{array}
\nonumber
\end{equation}
This immediately yields (\ref{1.5}).

For the proof of (\ref{1.6}) we recall the following inequality (see \cite{IS})
\begin{align*}
(\gamma-\gamma_0 I_3)\mbox{\boldmath $A$}\cdot\mbox{\boldmath $A$}
&+\gamma_0(\mbox{\boldmath $A$}-\mbox{\boldmath $B$})
\cdot (\mbox{\boldmath $A$}-\mbox{\boldmath $B$})
\\
&
=\gamma\mbox{\boldmath $A$}\cdot\mbox{\boldmath $A$}
-2\gamma_0\mbox{\boldmath $A$}\cdot\mbox{\boldmath $B$}
+\gamma_0\mbox{\boldmath $B$}\cdot\mbox{\boldmath $B$}
\\
&
=
\left\vert\gamma^{1/2}\mbox{\boldmath $A$}-\gamma_0\gamma^{-1/2}\mbox{\boldmath $B$}\right\vert^2
+\gamma_0\mbox{\boldmath $B$}\cdot\mbox{\boldmath $B$}
-\gamma_0^2\gamma^{-1/2}\mbox{\boldmath $B$}\cdot\gamma^{-1/2}\mbox{\boldmath $B$}
\\
&
\geq
\gamma_0\gamma\gamma^{-1/2}\mbox{\boldmath $B$}\cdot\gamma^{-1/2}\mbox{\boldmath $B$}
-\gamma_0^2\gamma^{-1/2}\mbox{\boldmath $B$}\cdot\gamma^{-1/2}\mbox{\boldmath $B$}
\\
&
=\gamma_0(\gamma-\gamma_0 I_3)\gamma^{-1/2}\mbox{\boldmath $B$}
\cdot\gamma^{-1/2}\mbox{\boldmath $B$},
\end{align*}
where $\mbox{\boldmath $A$}$ and $\mbox{\boldmath $B$}$ are real vectors.
Applying this to the first and second term 
in the right-hand side on (\ref{A.4}),  we obtain
\begin{equation}
\begin{array}{ll}
\displaystyle
-\int_{\Bbb R^3}f Rdx &
\displaystyle
\ge 
\int_{\Bbb R^3}\gamma_0(\gamma-\gamma_0I_3)\gamma^{-1/2}\nabla v\cdot\gamma^{-1/2}\nabla v dx
\\[5mm]
\displaystyle
&
\displaystyle
\,\,\,
+
e^{-\tau T}\left(\int_{\Bbb R^3}FRdx+\int_{\Bbb R^3}Fvdx\right).
\end{array}
\nonumber
\end{equation}
Now applying (\ref{A.5}), (\ref{A.8}) and (\ref{A.9}) to the second term in this
right-hand side we obtain (\ref{1.6}).

%
%
%

%
%
%
\end{document}